\documentclass[a4paper,12pt,reqno]{amsart}
\usepackage{amsmath,amsfonts}
\usepackage[numbers]{natbib}

\RequirePackage[mathscr]{eucal}

\DeclareMathOperator{\expect}{{\mathbb E}}
\DeclareMathOperator{\prob}{{\mathbb P}}

\DeclareMathOperator{\id}{Id}

\makeatletter
\newcommand{\eqcolon}{\mathrel{\mathord{=}\raise.2\p@\hbox{:}}}
\newcommand{\coloneq}{\mathrel{\raise.2\p@\hbox{:}\mathord{=}}}
\makeatother
\newcommand{\R}{\mathbb{R}}
\newcommand{\wt}{\widetilde}
\newcommand{\norm}[1]{\|#1\|}

\newtheorem{lemma}{Lemma}
\newtheorem{proposition}{Proposition}
\newtheorem{corollary}{Corollary}
\newtheorem{theorem}{Theorem}
\newtheorem{remark}{Remark}

\newcommand{\RR}{\mathbb{R}^3}
\newcommand{\CC}{\mathscr{C}}
\newcommand{\XT}{\mathbf{X}_T}
\newcommand{\XTO}{\mathbf{X}_{T_0}}
\newcommand{\DXT}{\mathbf{D}_{X,T}}
\newcommand{\DXoT}{\mathbf{D}_{X,\overline{T}}}

\newcommand{\DX}{\mathscr{D}_{X}}
\newcommand{\DWTX}{\mathscr{D}_{\wt X}}
\newcommand{\DXTO}{\mathbf{D}_{X,T_0}}
\newcommand{\DWTXTO}{\mathbf{D}_{\wt X,T_0}}

\begin{document}
\title[Evolution of a Random Vortex]{The evolution of a random  vortex filament}
\date{July 2004}
\author{Hakima Bessaih \and Massimiliano Gubinelli \and Francesco Russo}
\address{
Dept. of Mathematics\\
University of Wyoming\\
Ross Hall 210 --
Laramie WY 82071, USA}
\email{Bessaih@uwyo.edu}
\address{Institue Galil\'ee, Mathematiques\\
Universit\'e Paris 13\\
99, av. J.-B. Cl\'ement - F-93430 Villetaneuse, FRANCE 
}
\email{gubinell@math.univ-paris13.fr}
\address{Institue Galil\'ee, Mathematiques\\
Universit\'e Paris 13\\
99, av. J.-B. Cl\'ement - F-93430 Villetaneuse, FRANCE }
\email{russo@math.univ-paris13.fr}

\abstract
We study an evolution problem in the space of continuous loops in a
three-dimensional Euclidean space modelled upon the dynamics of vortex
lines in 3d incompressible and inviscid fluids. We establish existence
of a local solution starting from H\"older regular loops with index
greater than $1/3$.  When the H\"older regularity of the initial
condition $X$ is smaller or equal $1/2$ we require $X$ to be a
\emph{rough path} in the sense of Lyons~\cite{Lyons,LyonsBook}. The solution
will then live in an appropriate space of rough paths. In particular
we can construct (local) solution starting from almost every Brownian loop.
\endabstract
\maketitle

\section{Introduction}

The aim of this work is to study the well-posedness of the evolution
problem for a model of a random \emph{vortex filament} in
three dimensional incompressible fluid. If $u$ is the velocity field
of the fluid, the vorticity $\omega : \mathbb{R}^3 \to \mathbb{R}^3$
is a solenoidal field defined as $\omega = \text{curl}\, u$. 
A vortex filament is a field of vorticity $\omega$  
which is strongly concentrated around a three-dimensional closed curve
$\gamma$ described parametrically as a continuous function $\gamma :
[0,1] \to \mathbb{R}^3$ such that $\gamma_0 = \gamma_1$. Ideally, neglecting the transverse size of
the filament, we can describe the vorticity field $\omega^\gamma$
generated by $\gamma$ formally as the distribution
\begin{equation}\label{formal}
  \omega^\gamma(x)=\Gamma\int_0^1\delta(x-\gamma_\xi) d_\xi \gamma_\xi,
  \qquad x \in\mathbb{R}^3
\end{equation}
where $\Gamma > 0$ is the intensity of vorticity. In
$\mathbb{R}^3$, the velocity field associated to $\omega$ can be
reconstructed with the aid of the Biot-Savart formula:
\begin{equation}\label{biot-savart}
u^\gamma(x)=-\frac{1}{4\pi} \int_{\R^3}\frac{(x-y)}{|x-y|^3}\wedge \omega^\gamma( y)\,d y.
\end{equation}
which is the solution of $\text{curl}\, u^\gamma = \omega^\gamma$ with
enough decay at infinity.

Then
\begin{equation}\label{biot-savart2}
u^\gamma(x)=-\frac{\Gamma}{4\pi}
\int_0^1\frac{(x-\gamma_\xi)}{|x-\gamma_\xi|^3}\wedge \, d_\xi \gamma_\xi.
\end{equation}
where $a \wedge b$ is the vector product of the vectors $a,b \in
\mathbb{R}^3$. The evolution in time of the infinitely thin vortex filament is
obtained by imposing that the curve $\gamma$ is transported by the
velocity field $u^\gamma$:
\begin{equation}
  \label{eq:evolution1}
 \frac{d}{dt} \gamma(t)_\xi = u^{\gamma(t)}(\gamma(t)_\xi)
, \qquad \xi
 \in [0,1]  
\end{equation}
and this gives the initial value problem
\begin{equation}
  \label{eq:evolution2}
 \frac{d}{dt} \gamma(t)_\xi 
= -\frac{\Gamma}{4\pi}
\int_0^1\frac{(\gamma(t)_\xi-\gamma(t)_\eta)}{|\gamma(t)_\xi-\gamma(t)_\eta|^3}\wedge \, d_\eta \gamma(t)_\eta
, \qquad \xi
 \in [0,1]  
\end{equation}
Even if the curve $\gamma$ is smooth this expression is not well
defined since the integral is divergent if $\gamma$ has non-zero
curvature. 
 
To overcome this divergence a natural approach is that of
Rosenhead~\cite{Rosenhead}, who suggested the following approximate equation of
motion based on a regularized kernel
\begin{equation}
  \label{eq:evolution-reg}
 \frac{d}{dt} \gamma(t)_\xi 
= -\frac{\Gamma}{4\pi}
\int_0^1\frac{(\gamma(t)_\xi-\gamma(t)_\eta)}{[|\gamma(t)_\xi-\gamma(t)_\eta|^2
  + \mu^2]^{3/2}}\wedge \, d_\eta \gamma(t)_\eta
, \qquad \xi
 \in [0,1]  
\end{equation}
for some $\mu>0.$
This model has clear advantages and been used in some  numerical
calculation of  aircraft trailing vortices by Moore~\cite{Moore}.

We will consider a generalization of the Rosenhead model 
where the function $\gamma$ is not necessarily smooth. This
is natural when we want to study models of
\emph{random} vortex filaments. Indeed for simple models of random
vortex lines the curve $\gamma$ is rarely smooth or even of bounded
variation. Here we imagine to take as initial condition a
typical trajectory of a Brownian loop (since the path must be closed)
or other simple models like fractional Brownian loops (to be described
precisely in Sec.~\ref{sec:fbl}). 
As we will see later, a major  problem is then the interpretation of 
equation~(\ref{eq:evolution-reg}).

The study of the dynamics of random vortex lines is suggested by some
work of A.~Chorin~\cite{Chorin} and G.~Gallavotti~\cite{Gallavotti}.
The main justification for the adoption of a probabilistic point of
view comes from two different directions. Chorin builds discrete models of
random vortex filaments to explain the phenomenology of turbulence by
the statistical mechanics of these coherent structures. Gallavotti
instead suggested the use of very irregular random functions to
provide a natural regularization of eq.~(\ref{eq:evolution2}). Both
approaches are on a physical level of rigor. 

On the mathematical side in the recent year there have been some
interest in the study of the statistical mechanics of continuous
models for vortex filaments. 
P.L.~Lions and A.J.~Majda~\cite{LionsMajda} proposed a
statistical model of quasi-3d random vortex lines which are constrained
to remain parallel to a given direction and thus cannot fold.
Flandoli~\cite{Flandoli} rigorously studied the problem of the
definition of the energy for a random vortex filament modeled over a
Brownian motion and Flandoli and Gubinelli~\cite{FlaGub} introduced a
probability measure over Brownian paths to study the statistical
mechanics of vortex filaments. The study of the energy of filament
configurations has been extended to models based on fractional
Brownian motion by Flandoli
and Minnelli~\cite{FlaMin} and Nualart, Rovira and Tindel~\cite{NRT}.
Moreover a model of Brownian vortex
filaments capable of reproducing the multi-fractal character~\cite{Frisch} of
turbulent velocity fields has been introduced in~\cite{FlaGub2}.

The problem~(\ref{eq:evolution-reg}) define a natural flow on three
dimensional closed curves. The study of this kind of flows has been
recently emphasized by Lyons and Li in~\cite{LiLyons} where they study
a class of flows of the form
\begin{equation}
  \label{eq:driver-flow}
\frac{dY}{dt} = F(I_f(Y)), \qquad Y_0 = X  
\end{equation}
where $Y$ takes values in a (Banach) space of functions, $F$ is a
smooth function and $I_f(Y)$ is an It\^o map, i.e. the map $Y
\mapsto Z$ where $Z$ is the solution of
the differential equation
$$
dZ_\sigma = f(Z_\sigma) dY_\sigma
$$
driven by the path $Y$. They prove that, under suitable conditions,
$I_f$ is a smooth map and then that eq.~(\ref{eq:driver-flow})
has (local) solutions and thus as a by-product that $F \circ I_f$ can
be effectively considered a vector-field on a space of paths. 

Our evolution
equation does not match the structure of the flows considered by Lyons
and Li. A very important difference is that eq.~(\ref{eq:driver-flow}),
under suitable assumptions on the initial condition $X$ (e.g.  $X$ a
semi-martingale path), can be solved with standard tools of
stochastic analysis (essentially It\^o stochastic calculus) while
eq.~(\ref{eq:evolution-reg}) has a structure which is not well adapted
to a filtered probability space and prevents even to (easily) set-up
the problem in a space of semi-martingale paths. To our opinion this
peculiarity makes the problem interesting from the point of view of
stochastic analysis and was one of our main motivation to start its
study. 
Using Lyons' rough paths we will show that it is possible to give a
meaning to the evolution problem~(\ref{eq:evolution-reg}) starting form a (fractional or
standard) Browian loop and that this problem has always a local
solution (recall that the existence of a global solution, to our
knowledge, has not been proven even in the case of a smooth curve). 

The paper is organized as follows: in Sec.~\ref{sec:model} we describe
precisely the model we are going to analyze and we make some
preliminary observations on the structure of covariation  of the
solution (in the
sense of stochastic analysis and assuming the initial condition has
finite covariations). 
Next, we introduce the functional spaces in which we will set-up
the problem of existence of solutions. In Sec.~\ref{sec:young} we
build a local solution for initial conditions which are H\"older
continuous with exponent greater that $1/2$ and for which the line
integrals can be understood \emph{\`a la} Young~\cite{young}. In
Sec.~\ref{sec:rough} we build a solution in a class of rough paths
(introduced in~\cite{Gubinelli}) for initial conditions which are
rough paths of H\"older regularity greater that $1/3$ (which
essentially are $p$-rough paths for $p < 3$, in the terminology
of~\cite{LyonsBook}). Moreover we prove that the solution is Lipshitz
continuous with respect to the initial data. 
Finally, in Sec.~\ref{sec:random} we apply these
results to obtain the evolution of random initial
conditions of Brownian loop type or its fractional
variant. 
Appendix~\ref{sec:app-proofs} collect the proof of some lemmas.

\section{The Model}
\label{sec:model}
\subsection{The evolution  equation}
Our aim is to start a study of the  tridimensional evolution of
random vortex filaments by the analysis of the well-posedness of the
regularized dynamical equations.
Inspired by the Rosenhead model~(\ref{eq:evolution-reg}) we will be
 interested in studying the  evolution described by
\begin{equation}
  \label{eq:evolution}
\frac{\partial Y(t)_\xi}{\partial t} = V^{Y(t)}(Y(t)_\xi), \qquad Y(0)
= X
\end{equation}
with initial condition $X$ belonging to the set $\CC$ of  closed and
continuous curve in $\RR$ parametrized by $\xi \in [0,1]$. For any $Z
\in \CC$,  $V^{Z}$ is the vector-field given by the line integral
\begin{equation}
  \label{eq:evolutionV}
  V^{Z}(x) = \int_Z A(x-y) dy = \int_0^1 A(x-Z_\xi) dZ_\xi, \qquad x \in \mathbb{R}^3
\end{equation}
where $A : \RR \to \RR \otimes \RR$ is a matrix-valued field.
In this setting the Rosenhead model is obtained by taking $A$ of the form
$$
A(x)^{ij} = -\frac{\Gamma}{4\pi} \sum_{k=1,2,3} \epsilon_{ijk}
\frac{x^k}{[|x|^2+\mu^2]^{3/2}} 
,\qquad i,j = 1,2,3, \quad x \in \mathbb{R}^3
$$
where $\epsilon_{ijk}$ is the completely antisymmetric tensor in
$\mathbb{R}^3$ normalized such that $\epsilon_{123}=1$, $\mu > 0$ is a
fixed constant and $(x^k)_{k=1,2,3}$ are the components of the vector $x
\in \mathbb{R}^3$.

\subsection{A first approach using covariations for random initial conditions}
\label{sec:quad-var}
Even if the kernel of the paper is fully pathwise,  
before studying existence and uniqueness problem (in some sense to be
precised),  we would like to insert a preamble
concerning the stability of  the quadratic variation (with respect to
the parameter) of a large class of solutions.

Let $(\Omega, {\mathcal F}, \mathrm{P})$ be a probability space.
In order to simplify a bit the proofs we have chosen to use the notion
of covariation introduced, for instance, in~\cite{rv95}.  

 Given two processes $X = (X_\xi)_{\xi \in [0,1]} $ and 
$Y = (Y_\xi)_{\xi \in [0,1]} $, the covariation  $[X,Y]$ is defined (if it exists) as the limit
of the sequence of functions 
$$ \xi \mapsto \int_0^\xi (X_{\rho+\varepsilon}-X_\rho) (Y_{\rho+\varepsilon}-Y_\rho)\frac{d\rho}{\varepsilon} $$
in the uniform convergence in probability sense (ucp).
If $X$ and $Y$ are classical continuous semi-martingales, 
 it is well-known that previous $[X,Y]$ coincides with the classical covariation.

 A vector $(X^1, \ldots, X^n)$ of stochastic processes
is said to have {\it all its mutual covariations} 
if $[X^i,X^j]$ exist for every $i,j = 1,\ldots, n$. 
Generally here we will consider $n = 3$. Moreover, given a matrix or vector $v$ 
we denote $v^*$ its transpose. 

It is sometimes practical to have a matrix notations.
If $M_\xi = \{m^{ij}_\xi\}_{i,j}$, 
$N_\xi = \{n^{ij}_\xi\}_{i,j}$, are matrices of stochastic
processes such they are compatible for the matrix
product, we denote
$$[M,N]_\xi = \left\{ \sum_{k=1}^n [m^{ik}, n^{kj}]_\xi\right\}_{i,j} $$

\begin{remark} \label{r1}
 The following result can be easily deduced
 from~\cite{rv95}.
Let $\Phi_1, \Phi_2$ be of class $C^1(\R^3; \R^3)$,
$X = (X^1, X^2, X^3)^*, Z= (Z^1,Z^2, Z^3)^* $
(understood as row vectors in the matrix calculus)
such that $(X,Z)$ has all its mutual covariations.
Then $(\Phi_1(X), \Phi_2(Z))$ has all its mutual covariations
and 
$$
[\Phi_1(X),\Phi_2(Z)^*]_\xi =  \int_0^\xi   (\nabla
\Phi_1)(Z_\rho) d[X,Z^*]_\rho (\nabla \Phi_2)(X_\rho)^*. 
$$
\end{remark}

\begin{remark} \label{r2}
In reality, we could have chosen the modified F\"ollmer~\cite{fo} approach
appearing in~\cite{erv}, based on discretization procedures for which the common reader
would be more accustomed.

In that  case the same results stated in Remark~\ref{r1} and
Proposition~\ref{PQVar} will be valid also in this discretization framework.
We recall briefly that context.

Consider a family of subdivisions
$ 0 = \xi^n_0 <  \ldots < \xi^n_n =  1 $ of the interval $[0,1]$.
We will say that the mesh of the subdivision converges to  zero
if $|\xi^n_{i+1} -  \xi^n_{i}| $ go to zero as $n \to \infty$.

 In this framework,  the covariation  $[X,Y]$ is defined (if it exists) as the limit
of 
$$ \xi \mapsto \sum_{i=0}^{n-1}  (X_{\xi^n_{i+1} \wedge \xi }- X_{\xi^n_{i} \wedge \xi}) (Y_{\xi^n_{i+1} \wedge \xi}- Y_{\xi^n_{i} \wedge \xi}) $$
 in the ucp (uniform convergence in probability) 
sense with respect to $\xi$ 
and the limit does not depend on the chosen family of subdivisions. 
\end{remark}

\medskip
 
Suppose there exists a sub Banach space $B$ of  $\CC$ 
and let $V : (\gamma, y) \rightarrow V^\gamma (y)$
 be  a Borel map 
from $B  \times  \R^3$ to $\R^3$ such that 
\begin{itemize} \label{Ass}
\item[V1)] for fixed $ \gamma \in B  $,
$ y \mapsto V^\gamma (y)$   is $C^1_{\rm b} (\R^3; \R^3 )$;
\item[V2)] the application
$\gamma  \mapsto \|\nabla V ^\gamma\|_\infty$   is locally bounded from
$B$ to $\R$.
\end{itemize} 

The main motivation for this abstract framework comes from the setting
described in the following section. Indeed, as we will see, there
exists natural Banach sub-spaces $B$ of $\CC$ such that the map
$V$ defined as
$$ V^\gamma(x) = \int_0^1 A(x - \gamma_\xi) d^* \gamma_\xi,  $$
where $d^*$ denotes some kind of path integration defined for every
$\gamma \in B$, satisfy the above hypoteses V1) and V2). 

\begin{proposition} \label{PQVar} 
Suppose that a random field $(Y(t)_\xi)_{\xi \in [0,1], t \in [0,T]}$ with values in $B$ is a continuous solution of
\begin{equation} \label{EQVar}
 Y(t)_\xi = X_\xi + \int_0^t V^{Y(s)} (Y(s)_\xi) ds;\qquad \xi \in
 [0,1], t \in [0,T] 
\end{equation}
 with an
  initial condition $X$ having all its mutual covariations.
Then at each time $t \in [0,T]$ the path $Y(t)$ has all its mutual covariations.
Moreover
\begin{equation}
  \label{eq:EQVar2}
[Y(t), Y^*(t)]_\xi  =   \int_0^\xi  M(t)_\rho  d[X,X^*]_\rho (M(t)_\rho)^*  
\end{equation}
where
\begin{equation}
  \label{eq:Mprocess}
M(t)_\xi := \exp \left[\int_0^t
(\nabla V^{Y(s)})(Y(s)_\xi) ds \right].  
\end{equation}
\end{proposition}

\begin{remark} \label{RQVar} 
Note the following:
\begin{enumerate}
\item It is possible to adapt this proof to the situation where the
  solution exists up to a random time.
\item Since we are in the multidimensional   case, 
we recall that $E_{\xi}(t) = \exp \left[\int_0^t
Q_\xi(s) ds \right]$ is defined as the
matrix-valued function satisfying the differential equation
$$
E_\xi(t) \in \RR \otimes \RR , \quad \frac{d E_\xi(t)}{dt} = Q_\xi(t) E_\xi(t), \qquad
E_\xi(0) = \text{Id}.
$$
\item A typical case of initial condition of process having
all its mutual covariation is a 3-dimensional Brownian loop. 
\end{enumerate}
\end{remark}

\begin{proof}

For simplicity, we prolongate the processes $X$ parametrized by
$[0,1]$ setting  $X_\xi = X_1$, $\xi\ge 1$. 
Let $\varepsilon > 0$. 
For $t \in [0,T], \xi \in [0,1]$,
write 
$$ Z^\varepsilon (t)_\xi =  Y(t)_{\xi + \varepsilon}  -  Y(t)_\xi , \quad X^\varepsilon_\xi  =  X_{\xi + \varepsilon}   - X_\xi, $$
then
\begin{eqnarray*}
   Z^\varepsilon (t)_\xi &=&  X^\varepsilon_\xi + \int_0^t  [
  V^{Y(s)}  (Y(s)_{\xi + \varepsilon} )-  V^{Y(s)} (Y(s)_\xi )  ]ds         \\  
&=&     X^\varepsilon_\xi + 
      \int_0^t  ( \nabla V^{Y(s)})  (Y(s)_\xi ) Z^\varepsilon (s)_\xi
  \,ds\\
& &  +  \int_0^t R^\varepsilon (s)_\xi Z(s)_\xi\,ds
 \end{eqnarray*}
where 
$$ R^\varepsilon (s)_\xi = \int_0^t ds  \int_0^1 da \left [ (\nabla V^{Y(s)}) (Y(s)_\xi  + a Z^\varepsilon (s)_\xi) - 
 ( \nabla V^{Y(s)})  (Y(s)_\xi)  \right] $$
so that 
$$ \sup_{s \le T} \vert R^\varepsilon (s)_\xi \vert  \rightarrow 0, \quad {\rm }a.s.  $$
Therefore,
$$ Z^\varepsilon (t)_\xi  = \exp \left[ \int_0^t ds  (\nabla V^{Y(s)})
  (Y(s)_\xi)+  \int_0^t ds R^\varepsilon (s)_\xi \right]
 X^\varepsilon_\xi  \eqcolon M(t)^\varepsilon_\xi X^\varepsilon_\xi. $$

Multipling both sides by their transposed, dividing by $\varepsilon$ and integrating from $0$ to $y$
we get
$$ \int_0^\xi \frac{Z^\varepsilon(t)_\rho
  (Z^\varepsilon)^*(t)_\rho}{\varepsilon} d\rho = \int_0^\xi
  M(t)^\varepsilon_\rho 
 \frac{X^\varepsilon_\rho (X^\varepsilon)^*_\rho}{\varepsilon}
  (M(t)^\varepsilon_\rho)^* d\rho. $$

Then since, as $\varepsilon \to 0$, 
$$
M(t)^\varepsilon_\xi \to  \exp\left[\int_0^t (\nabla V^{Y(s)}(Y(s)_\xi)
  ds \right] = M(t)_\xi
$$
uniformly in $t$ and $\xi$ almost surely,  using Lebesgue dominated convergence theorem,
and similar arguments to Proposition~2.1 of~\cite{rv95},
it is enough to show that
$$ \xi \mapsto \int_0^\xi
\exp \left[ \int_0^t ds (\nabla V^{Y(s)}) (Y(s)_\rho ) \right]
   \frac{X^\varepsilon_\rho (X^\varepsilon_\rho)^*}{\varepsilon} \exp
   \left[ \int_0^t ds (\nabla V^{Y(s)}) (Y(s)_\rho ) \right]^* d\rho $$
converges ucp  to the right member of~(\ref{eq:EQVar2}). 

This is obvious since 
$$ \int_0^\xi \frac{X^\varepsilon_\rho (X^\varepsilon_\rho)^*}{\varepsilon} d\xi \rightarrow  [X,X^*]_\xi$$
ucp with respect to $\xi \in [0,1]$ and so, modulo extraction of a  subsequence, we can
make use of the weak $\star$-topology. 
\end{proof}

\begin{remark}\label{Ncov}
Suppose that  the initial condition  has all its $n$-mutual covariations
 $n\ge 3$, see for this~\cite{er}.
Proceeding in a similar way as above, it is possible to show 
that $Y(t)$ has all its  mutual $n$-covariations.


A typical example of process having a strong finite 
$n$-variation is fractional Brownian motion with Hurst
index $H= 1/n$.
\end{remark}


\subsection{The functional space framework}

The filament evolution  problem when $X$ has finite-length has been previously studied
in~\cite{BerselliBessaih} where it is proved that under some regularity conditions
on $A$ there exists a unique local solution living in the space $H^1_c$ of
closed curves with $L^1$ derivative. 

We would like to be able to solve the Cauchy
problem~(\ref{eq:evolution}) starting from a  random curve $X$ like a
3d Brownian loop (since it must be closed) or a fractional Brownian
loop. In these cases $X$ is a.s. not in $H^1_c$ and, as a consequence, we need
 a sensible definitions to the path-integral appearing in
eq.~(\ref{eq:evolutionV}).

Even if $X$ is a Brownian loop do not exists a simple strategy to give
a well defined meaning to the evolution problem~(\ref{eq:evolution})
through the techniques of stochastic calculus. Indeed we could try to
define the integral in $V^Y$ as an It\^o or Stratonovich integral
which requires $Y$ to be a semi-martingale with respect to some
filtration $\mathcal{F}$ (e.g. the filtration generated by $X$). However we readily note that the problem has
no relationship with any natural filtration $\mathcal{F}$ since for example to
compute the velocity field $V^Y(x)$ in some point $x$ we need information about the whole trajectory of $Y$. 

A viable (and relatively straightforward) strategy is then to give a well defined
meaning to the problem  using a path-wise approach. 

We require that the initial data has
$\gamma$-H\"older regularity. 
When $\gamma >
1/2$ the line integral appearing in the
definition~(\ref{eq:evolutionV}) of the instantaneous velocity field
$V^{Y(t)}$ can be understood \emph{\`a la} Young~\cite{young}. The
corresponding results will be presented in
Sec.~\ref{sec:young}. 

When $1/2 \ge  \gamma > 1/3$ an appropriate
notion of line integral has been formulated by Lyons
in~\cite{Lyons,LyonsBook}. In Sec.~\ref{sec:rough} we will show that
given an initial $\gamma$-H\"older path $X$ (and
an associated \emph{area process} $\mathbb{X}^2$) there exists a  unique local
solution of the problem~(\ref{eq:integrated}) in the class $\DX$ of paths
\emph{weakly-controlled} by $X$. The class $\DX$ has been introduced
in~\cite{Gubinelli} to provide an alternative formulation of Lyons' theory of
integration and corresponds to paths $Z \in \CC$ which locally behaves as
$X$ in the sense that
$$
Z_\xi - Z_\eta = F_\eta (X_\xi-X_\eta) + \text{O}(|\xi-\eta|^{2\gamma})
$$
where $F \in C([0,1],\RR \otimes \RR)$ is a path taking values in the
bounded endomorphisms of $\RR$.

In particular these results provide solutions of the problem when $X$
is a fractional Brownian loop of Hurst-index $H > 1/3$ (see Sec.\ref{sec:random}).

\section{Evolution for $\gamma$-H\"older curves with $\gamma > 1/2$}
\label{sec:young}

\subsection{Setting and notations}

For any $X \in \CC$ let
\begin{equation*}
  \|X\|_{\gamma} :=
\sup_{\xi,\eta\in[0,1]} \frac{|X_\xi-X_\eta|}{|\xi-\eta|^\gamma},
\qquad   \|X\|_{\infty} :=
\sup_{\xi \in[0,1]} |X_\xi|
\end{equation*}
and
\begin{equation*}
  \|X\|_{\gamma}^{*} := \norm{X}_{\infty}+
  \norm{X}_{\gamma}
\end{equation*}
Denote $\CC^\gamma$ the set of paths $X \in \CC$ such that
$\norm{X}_{\gamma}^*< \infty$.

All along this section we will assume that $\gamma$ is a fixed number
greater than $1/2$.
In this case the following result states that there exists a unique
extension to the Riemann-Stieltjes integral $\int f dg$ defined for
smooth functions $f,g$ to all  $f,g \in \CC^\gamma$.

\begin{proposition}[Young's integral]
\label{lemma:young}
Let $X,Y \in \CC^\gamma$, then
$
   \int_\eta^\xi X_\rho d Y_\rho
$
is well defined, coincide with the Riemann-Stieltjes integral when
the latter exists and satisfy the following bound
\begin{equation*}
  \left|\int_\eta^\xi (X_\rho-X_\eta) dY_\rho \right| \le C_{\gamma} \|X\|_{\gamma} \|Y\|_{\gamma} |\xi-\eta|^{2\gamma}
\end{equation*}
for all $\xi,\eta \in [0,1]$ where $C_\gamma \ge 1$ is a constant depending
only on $\gamma$.
\end{proposition}
\begin{proof}
See e.g.\cite{young,Lyons}.
\end{proof}

It will be convenient to introduce the integrated form
of~(\ref{eq:evolution}) as
\begin{equation}
  \label{eq:integrated}
 Y(t)_\xi = X_\xi + \int_0^t V^{Y(s)}(Y(s)_\xi) ds
\end{equation}

Consider the Banach space $\XT \coloneq C([0,T],\CC^\gamma)$ with norm
$$
\|Y\|_{\XT} \coloneq \sup_{t \in [0,T]} \|Y(t)\|_{\gamma}^*,\qquad Y
\in \XT.
$$

Solutions of~(\ref{eq:integrated}) will then be found as fixed points
of the non-linear map $F : \XT \to \XT$ defined as
\begin{equation}
  \label{eq:fmap}
F(Y)(t)_\xi \coloneq X_\xi + \int_0^t V^{Y(s)}(Y(s)_{\xi}) \, ds  ,\qquad t
\in[0,T], \xi \in [0,1].
\end{equation}
where the application $Z \mapsto V^Z$ is defined for any $Z \in
\CC^\gamma$ as in eq.~(\ref{eq:evolutionV}) with the line integral
understood according to proposition~\ref{lemma:young} and with the
matrix field $A$ satisfying regularity conditions which will be
shortly specified.

On $m$-tensor field
 $\varphi: \RR \to (\RR)^{\otimes m}$ and for any integer $n \ge 0$ we define the following
norm:
$$
\norm{\varphi}_n \coloneq \sum_{k=0}^n \norm{\nabla^k\varphi }
\qquad \text{where}
\qquad
\norm{\varphi} \coloneq \sup_{x \in \RR} |\varphi(x)|.
$$
where the norm $|M|$ of a matrix $M \in \RR \otimes \RR$ or more
generally of a $n$-tensor  $M \in (\RR)^{\otimes n}$ is given by
$$
|M| = \sum_{i_1=1}^3 \cdots \sum_{i_n=1}^3 |M^{i_1\cdots i_n}|
$$
with $(M^{i_1\cdots i_n})_{i_1,\dots,i_n}$ the components of the tensor in the canonical
basis of $\RR$.

Then we can state:
\begin{theorem}
\label{th:sol-young}
Assume $\|\nabla A\|_2 < \infty$ and $X \in \CC^\gamma$. Then there
exists a time $T_0$ depending only on $\|\nabla A\|_1, \|X\|_\gamma^*,
\gamma$ such that the equation~(\ref{eq:integrated}) has a unique
solution bounded in $\CC^\gamma$.
\end{theorem}
\begin{proof} Consider the initial condition $X$ fixed.
According to 
lemma~\ref{lemma:existence-young} below, there exists $T_0 >0 $ and $B_{T_0}>0$ depending
on  $\|\nabla A\|_1$, $\|X\|_\gamma^*$,
$\gamma$ such that the set $\{Y \in \XTO : Y(0) = X, \|Y\|_{\XTO} \le B_{T_0}\}$ is invariant under $F$. 
Lemma~\ref{lemma:uniqueness-young} then assert that,  provided $\|\nabla A\|_2 < \infty$, the map $F$ is a
strict contraction over a smaller time interval $[0,\overline{T}]$
with $\overline{T}
\le T_0$. Proceding by induction on the intervals $[0,\overline{T}]$,
$[\overline{T},2\overline{T}]$, etc\dots it is possible to construct
the unique solution of the evolution problem
 up to the time $T_0$. 
\end{proof}

Before giving the lemmas used in the proof we will state two useful
results. The first is just a straighforward computation on H\"older functions, the second will allow to control the velocity field $V^Y$ in terms of
the regularity of $Y$ and of $A$.

\begin{lemma}
\label{lemma:stability-young}
Let $Y, \wt Y \in \CC^\gamma$ and $\varphi \in C^2$, then
\begin{equation}
  \label{eq:young-varphi-y-1}
 \norm{\varphi(Y)}_\gamma \le \norm{\nabla \varphi} \norm{Y}_\gamma 
\end{equation}
and
\begin{equation}
  \label{eq:young-varphi-y-2}
 \norm{\varphi(Y)-\varphi(\wt Y)}_\gamma \le \norm{\nabla
 \varphi}_1(1+ \norm{Y}_\gamma )\norm{Y-\wt Y}_\gamma
\end{equation}
\end{lemma}
\begin{proof} See Appendix, section~\ref{sec:proof-bounds-stab-young}.
\end{proof}

\begin{lemma}
\label{lemma:boundsyoung}
Let $Y,\wt Y \in \CC^\gamma$.
For any integer $n \ge 0$  the following estimates holds:
\begin{equation}
\label{eq:boundv}
  \begin{split}
\|\nabla^n V^{Y}\| \le C_\gamma \norm{\nabla^{n+1} A} \norm{Y}_{\gamma}^2
  \end{split}
\end{equation}
\begin{equation}
  \label{eq:boundvdiff}
\|\nabla^n V^{Y}- \nabla^n V^{\wt Y}\| \le
C_\gamma
\norm{\nabla^{n+1} A}_{1}(\norm{Y}_{\gamma}+\norm{\wt Y}_{\gamma}+\norm{\wt
  Y}_{\gamma} \norm{ Y}_{\gamma})  \norm{Y-\wt Y}_{\gamma}^*
\end{equation}
\end{lemma}
\begin{proof} See Appendix, section~\ref{sec:proof-bounds-young}.
\end{proof}

\subsection{Local existence and uniqueness}

\begin{lemma}
\label{lemma:existence-young}
Assume $\|\nabla A\|_1 < \infty$.
For any initial datum $X \in \CC^\gamma$
there exists a time $T_0>0$ such that for any
 time $T < T_0$ the set
$$
Q_T \coloneq \{Y \in \XT : Y(0) = X, \norm{Y}_{\XT} \le B_T \}
$$
where $B_T$ is a suitable constant, is invariant under $F$.
\end{lemma}
\begin{proof}
Let us compute
\begin{equation*}
  \begin{split}
    |F(Y)(t)_\xi| & \le |X_\xi| + \int_0^t |V^{Y(s)}(Y(s)_\xi)| ds
\le |X_\xi| + \int_0^t \norm{V^{Y(s)}}_\infty ds
  \end{split}
\end{equation*}
so
\begin{equation*}
  \begin{split}
  \|F(Y)(t)\|_{\infty}&  \le \|X\|_{\infty} +  \int_0^t
   \|V^{Y(s)}\| \, ds
\\& \le      \|X\|_{\infty} + C_\gamma \norm{\nabla A}\int_0^T  \|Y(s)\|_{\gamma}^2
 \, ds
\\& \le      \|X\|_{\infty} + T C_\gamma \norm{\nabla A}  \|Y\|_{\XT}^2
  \end{split}
\end{equation*}

The $\gamma$-H\"older norm of the path $F(Y)(t)$ can be estimated in a
similar fashion
\begin{equation*}
  \begin{split}
  \|F(Y)(t)\|_{\gamma} & \le \|X\|_{\gamma} +  \int_0^t
   \|V^{Y(s)}(Y(s)_\cdot)\|_\gamma \, ds
\\ & \le   \|X\|_{\gamma} +  \int_0^t
   \|\nabla V^{Y(s)}\| \|Y(s)\|_\gamma \, ds   
\\ & \le   \|X\|_{\gamma} +  C_\gamma  \norm{\nabla^2 A} \int_0^T \|Y(s)\|^3_\gamma
  \, ds   
\\ & \le   \|X\|_{\gamma} +  T C_\gamma  \norm{\nabla^2 A}  \|Y\|^3_{\XT}
  \end{split}
\end{equation*}
where we used eq.(\ref{eq:young-varphi-y-1}) in the second line and
eq.(\ref{eq:boundv}) in the third.

Then
\begin{equation*}
  \norm{F(Y)}_{\XT} \le \norm{X}_\gamma^* + C_\gamma T \norm{\nabla
  A}_1 \norm{Y}_{\XT}^2 (1+\norm{Y}_{\XT})
\end{equation*}
Let $B_T$ be a solution of
\begin{equation*}
B_T \le \norm{X}_\gamma^* + C_\gamma T \norm{\nabla A}_1 B_T^2 (1 + B_T)
\end{equation*}
which exists for any $T \le T_0$ where $T_0$ is a constant depending
only on $\norm{X}_\gamma^*$, $\norm{\nabla A}_1$ and $\gamma$.
Then if $\|Y\|_{\XT} \le B_T$  we have
$\|F(Y)\|_{\XT} \le B_T$ and $Q_{T}$ is invariant under $F$.
\end{proof}

Given another initial condition $\wt X \in \CC^\gamma$ consider the
associated map
\begin{equation}
  \label{eq:fmap-tilde}
\wt F(Y)(t)_\xi = \wt X_\xi + \int_0^t V^{Y(s)}(Y(s)_{\xi}) \, ds  ,\qquad t
\in[0,T],\, \xi \in [0,1].
\end{equation}

\begin{lemma}
\label{lemma:uniqueness-young}
 Assume $\|\nabla A\|_2 < \infty$. We have
 \begin{equation}
   \label{eq:lipshitz-bound-young}
\|F(Y)-\wt F(\wt Y)\|_{\XT}   \le \|X-\wt X\|_{\gamma}^* +
C_\gamma T \norm{\nabla A}_2  
 (1+\norm{Y}_{\XT}+\norm{\wt Y}_{\XT})^3
 \norm{Y-\wt
 Y}_{\XT}   
 \end{equation}
In particular, there exists a time $\overline{T} \le T_0$ such that the
map $F$ is a strict contraction on $Q_{\overline{T}}$.
\end{lemma}
\begin{proof}
We proceed as follows: take $Y,\wt Y \in \CC^\gamma$, then
\begin{equation*}
  \begin{split}
& \norm{F(Y)(t)-\wt F(\wt Y)(t)}_\infty \le \norm{X-\wt X}_\infty + \int_0^T\left[
\norm{V^{Y(s)}(Y(s)_\cdot) - V^{Y(s)}(\wt Y(s)_\cdot)}_\infty \right.
\\ & \left. \qquad \qquad \qquad  \qquad \qquad  \qquad \qquad +
\norm{V^{Y(s)}(\wt Y(s)_\cdot) - V^{\wt Y(s)}(\wt Y(s)_\cdot)}_\infty    
\right]\, ds
\\  & \qquad\le\norm{X-\wt X}_\infty + \int_0^T \left[ \norm{\nabla V^{Y(s)}} \norm{Y(s)-\wt
    Y(s)}_\infty
+ 
\norm{V^{Y(s)} - V^{\wt Y(s)}}    
\right]\, ds
\\  & \qquad\le \norm{X-\wt X}_\infty + C_\gamma \int_0^T \left[ \norm{\nabla^2
    A} \norm{Y(s)}_\gamma^2 \norm{Y(s)-\wt
    Y(s)}_\infty
\right.
\\ & \left.\qquad \qquad + 
\norm{\nabla A}_1 (\norm{Y(s)}_\gamma+\norm{\wt Y(s)}_\gamma+\norm{\wt
  Y(s)}_\gamma\norm{Y(s)}_\gamma)    \norm{Y(s)-\wt Y(s)}_\gamma^*
\right]\, ds
\\  & \qquad\le \norm{X-\wt X}_\infty +C_\gamma T \norm{\nabla A}_1 \norm{Y-\wt Y}_{\XT}
 (1+\norm{Y}_{\XT}+\norm{\wt Y}_{\XT})^2
  \end{split}
\end{equation*}
and
\begin{equation*}
  \begin{split}
& \norm{F(Y)(t)-\wt F(\wt Y)(t)}_\gamma \le \norm{X-\wt X}_\gamma + \int_0^T\left[
\norm{V^{Y(s)}(Y(s)_\cdot) - V^{Y(s)}(\wt Y(s)_\cdot)}_\gamma \right.
\\ & \left. \qquad \qquad \qquad  \qquad \qquad  \qquad \qquad +
\norm{V^{Y(s)}(\wt Y(s)_\cdot) - V^{\wt Y(s)}(\wt Y(s)_\cdot)}_\gamma    
\right]\, ds
\\  & \qquad\le  \norm{X-\wt X}_\gamma +  \int_0^T \left[ \norm{\nabla V^{Y(s)}}_1 \norm{Y(s)-\wt
    Y(s)}_\gamma^*
(1+\norm{Y(s)}_\gamma)
\right.
\\ & \left.\qquad \qquad \qquad
+ 
\norm{\nabla V^{Y(s)} - \nabla V^{\wt Y(s)}}\norm{\wt Y(s)}_\gamma    
\right]\, ds
\\  & \qquad\le \norm{X-\wt X}_\gamma +  C_\gamma \int_0^T \left[ \norm{\nabla^2
    A}_1\norm{Y(s)}_\gamma^2 (1+\norm{Y(s)}_\gamma) \norm{Y(s)-\wt
    Y(s)}_\gamma^*
\right.
\\ & \left.\qquad \qquad + 
\norm{\nabla^2 A}_1 (\norm{Y(s)}_\gamma+\norm{\wt Y(s)}_\gamma+\norm{\wt
  Y(s)}_\gamma\norm{Y(s)}_\gamma) \norm{\wt Y(s)}_\gamma    \norm{Y(s)-\wt Y(s)}_\gamma^*
\right]\, ds
\\  & \qquad\le  \norm{X-\wt X}_\gamma +  C_\gamma T \norm{\nabla^2 A}_1 \norm{Y-\wt Y}_{\XT}
 (1+\norm{Y}_{\XT}+\norm{\wt Y}_{\XT})^3
  \end{split}
\end{equation*}

Putting together the two estimates we get:
\begin{equation*}
  \begin{split}
\|F(Y)-F(\wt Y)\|_{\XT}  & \le  \norm{X-\wt X}^*_\gamma +
C_\gamma T \norm{\nabla A}_2  
 (1+\norm{Y}_{\XT}+\norm{\wt Y}_{\XT})^3
 \norm{Y-\wt
 Y}_{\XT}
  \end{split}
\end{equation*}
There exists $\overline{T} \le T_0$ such that
$$
C_\gamma \overline{T} \norm{\nabla A}_2 (1 + 2 B_{\overline{T}})^3 \eqcolon
\alpha < 1
$$
then if $Y, \wt Y \in Q_{\overline T}$ we have
\begin{equation*}
\|F(Y)-F(\wt Y)\|_{\XT}   \le  
C_\gamma T \norm{\nabla A}_2 
 (1+2 B_{\overline T})^3  \norm{Y-\wt Y}_{\XT} = \alpha  \norm{Y-\wt Y}_{\XT}
\end{equation*}
so that  $F$   is a strict contraction on $Q_{\overline{T}}$ with a unique fixed-point.
\end{proof}

\begin{remark}
Here and in the proofs for the case $\gamma > 1/3$ some conditions on
$A$ can be slightly relaxed using better estimates. For example, in
the proof of lemma~\ref{lemma:existence-young} the condition $\|\nabla
A\|_1 < \infty$ can be relaxed to require $\nabla A$ to be a H\"older
continuous function of index $(1-\gamma+\epsilon)/\gamma$ for some
$\epsilon > 0$, etc\dots\ However these refinements are not able to improve
qualitatively the results.
\end{remark}

\subsection{Dependence on the initial condition}

Denote with $T_0(\|X\|_\gamma^*)$ the existence time of the
solution build in Theorem~\ref{th:sol-young}, where we have considered
explicitly its dependence on the norm of the initial condition. For
any $r > 0$ let
$\mathscr{B}(\CC^\gamma;r)$ the open ball of $\CC^\gamma$ with radius $r$
and centered in zero.
Now, fix $r > 0$ and let $\Gamma : \mathscr{B}(\CC^\gamma;r) \to \XTO$,
where $T_0 = T_0(r)$, denote the solution of the evolution problem starting
from the initial condition $X \in \mathscr{B}(\CC^\gamma;r)$ and
living up to time $T_0(r)$.

\begin{theorem}
Under the conditions of Theorem~\ref{th:sol-young}, the map $X \mapsto \Gamma(X)$ is Lipshitz.  
\end{theorem}
\begin{proof}
Consider two initial conditions $X, \wt X \in
\mathscr{B}(\CC^\gamma;r)$. Note that $F(\Gamma(X)) = \Gamma(X)$ in
$\XTO$. By Lemma~\ref{lemma:uniqueness-young} we
have that, for $T < T_0$,
\begin{equation*}
  \begin{split}
   & \|F(\Gamma(X))-\wt F(\Gamma(\wt X))\|_{\XT}   =
    \|F(\Gamma(X))-\wt F(\Gamma(\wt X))\|_{\XT}  
\\ & \quad \le \|X-\wt X\|_{\gamma}^* +
C_\gamma T \norm{\nabla A}_2  
 (1+\norm{\Gamma(X)}_{\XT}+\norm{\Gamma(\wt X)}_{\XT})^3
 \norm{\Gamma(X)-\Gamma(\wt
 X)}_{\XT}   
  \end{split}
\end{equation*}
Since the norm of the initial condition is bounded by $r$ in
$\CC^\gamma$, by Lemma~\ref{lemma:existence-young} there exists a
constant $B_{T_0}(r)$ such that $\|\Gamma(X)\|_{\XTO} \le B_{T_0}(r)$
for any $X \in  \mathscr{B}(\CC^\gamma;r)$. Then for $T < T_0$
sufficiently small, we have
\begin{equation*}
  \begin{split}
\|\Gamma(X)-\Gamma(\wt X)\|_{\XT}   
 \le (1-\alpha)^{-1} \|X-\wt X\|_{\gamma}^* 
  \end{split}
\end{equation*}
where
$$
\alpha \coloneq C_\gamma T \norm{\nabla A}_2  
 (1+\norm{\Gamma(X)}_{\XT}+\norm{\Gamma(\wt X)}_{\XT})^3 < 1;
$$
which gives the Lipshitz continuity of the map on $\XT$. By an easy
induction argument it follows the Lipshitz continuity on all $\XTO$
(see e.g.~\cite{Gubinelli}).
\end{proof}

\subsection{Blow-up estimate}

From the previous results it is clear that if the norm
$\|Y(t)\|_{\gamma}$ of a solution $Y$ with initial condition $X \in \CC^\gamma$ is bounded by some number $M$  in some interval
$[0,\overline T]$, then the solution can be extended on a strictly
larger interval $[0,\overline T + \delta_M]$ with $\delta_M$ depending
only on $M$ (and on the data of the problem). This implies that the
only case in which we cannot find a global solution (for any positive
time) is when there is some time $\hat t_\gamma(X)$ such that $\lim_{t\to \hat t_\gamma(X)
  -} \|Y(t)\|_{\gamma} = +\infty$. This time is an \emph{epoch of
  irregularity} for the evolution in the class $\CC^\gamma$. Near this epoch we can establish a
lower bound for the norm $\|Y(t)\|_{\gamma}$.

\begin{proposition}
\label{prop:blowup}
Assume $\hat t_\gamma(X)> 0$ is the smallest epoch of irregularity for a solution $Y$ in the class $\CC^\gamma$. Then we have
\begin{equation}
  \label{eq:blowup}
  \|Y(t)\|_{\gamma} \ge \frac{C}{(\hat t_\gamma(X) - t)^{1/2}}
\end{equation}
for any $t \in [0,\hat t_\gamma(X))$.
\end{proposition}
\begin{proof}
\begin{equation}
\label{eq:V-estimate}
  \begin{split}
\|Y(t)\|_\gamma - \|Y(s)\|_\gamma & \le \int_s^t
\|V^{Y(u)}(Y(u)_\cdot)\|_\gamma du
\\ &      \le \int_s^t \|\nabla V^{Y(u)}\|_\infty
\|Y(u)\|_\gamma du
  \end{split}
\end{equation}
and using Lemma~\ref{lemma:boundsyoung} we have
\begin{equation*}
\|Y(t)\|_\gamma - \|Y(s)\|_\gamma
   \le C \int_s^t \|Y(u)\|^3_\gamma du
\end{equation*}
for some constant $C$ depending only on $A$ and $\gamma$, so that
\begin{equation*}
\frac{d}{dt} \|Y(t)\|_\gamma     \le  C\|Y(t)\|^3_\gamma
\end{equation*}
letting $y(t) = \|Y(t)\|_{\gamma}$ and integrating the
differential inequality between times $t>s$ we obtain
$$
\frac{1}{y(s)^2}-\frac{1}{y(t)^2} \le 2 C  (t-s)
$$
now, assume that there exists a time $\hat{t}_\gamma(X)$ such that $\lim_{t\to
  \hat t_\gamma(X) -} y(t) = +\infty$, then for any $s < \hat{t}_\gamma(X)$ we have the
following lower bound for the explosion of the $\CC^\gamma$ norm of $Y$:
$$
\|Y(s)\|_{\gamma} = y(s)^{1/2} \ge \frac{1}{(2C)^{1/2} (\hat t_\gamma(X) - s)^{1/2}}.
$$
\end{proof}

The estimate~(\ref{eq:V-estimate}) used in the previous proof implies
also that
\begin{equation}
  \label{eq:V-estimate2}
 z(t) \le z(0) + \int_0^t \|\nabla V^{Y(s)}\|_\infty z(s)ds
\end{equation}
where  $z(t) = \sup_{s\in[0,t]}\|Y(s)\|_\gamma$.
By Gronwall lemma
$$
z(t) \le z(0) \exp\left(\int_0^t \|\nabla V^{Y(s)}\|_\infty ds \right)
.
$$

This bound allows the continuation of any solution on the interval
$[0,t]$ if the integral $\int_0^t \|\nabla V^{Y(s)}\|_\infty ds$ is
finite. Then if  $\hat t_\gamma(X)$ is the first irregularity epoch
in the class $\CC^\gamma$ we must have that $\hat t_\gamma(X) =
\hat t(X) = \sup_{1/2 < \gamma \le 1} \hat t_\gamma(X)$ for any $1/2 < \gamma \le 1$. Indeed is easy to see that
for any $t < \hat t(X)$ there exists a finite constant
$M_t$ such that $\sup_{s \in [0,t]}\|\nabla V^{Y(s)}\|_\infty \le M_t$.

\begin{corollary}
Let $X \in \CC^{\gamma_*}$ with $\gamma_* > 1/2$, then
for any $1/2 < \gamma \le \gamma_*$ there exists a unique solution $Y^\gamma
\in C([0,\hat t_\gamma(X)),\CC^\gamma)$
with initial condition $X$. Moreover the first irregularity epoch $\hat t_\gamma(X)$ for the solution in
$\CC^\gamma$ does not depend on $\gamma \ge \gamma^*$.
\end{corollary}



\section{Evolution for $\gamma > 1/3$}
\label{sec:rough}

\subsection{Rough path-integrals}
When $\gamma \le 1/2$ there are difficulties in defining the
path-integral appearing in the expression~(\ref{eq:evolutionV}) for
the velocity field $V$. A successful
approach to such irregular integrals has been found by T.~Lyons to
consist in enriching the notion of \emph{path} (see
e.g.~\cite{LyonsBook,Lyons} and for some recent contributions~\cite{friz,Gubinelli,feyel}).

For any  $\gamma > 1/3$,
a $\gamma$-\emph{rough path}
(of degree two) is a couple $\mathbb{X} = (X,\mathbb{X}^2)$ where $X \in
\CC^\gamma$ and  $\mathbb{X}^2 \in
C([0,1]^2 , \RR\otimes \RR)$ is a matrix-valued function (called the \emph{area process}) on the square
$[0,1]^2$  verifying the following
compatibility condition with $X$:
\begin{equation}
  \label{eq:compatibility}
 \mathbb{X}^{2,ij} _{\xi\rho}- \mathbb{X}^{2,ij} _{\xi\eta}-
 \mathbb{X}^{2,ij} _{\eta\rho} = (X^i_\xi -X^i_\eta) (X^j_\eta-X^j_\rho),
 \qquad \xi,\eta,\rho \in [0,1]^{2}
\end{equation}
($i,j = 1,2,3$ are vector indexes)
and such that
\begin{equation}
  \label{eq:area-reg}
\norm{\mathbb{X}^2}_{2\gamma} \coloneq \sup_{\xi,\eta \in [0,1]} \frac{|\mathbb{X}^2_{\xi\eta}|}{|\xi-\eta|^{2\gamma}}< \infty.
\end{equation}

\begin{remark}
If $\gamma > 1/2$ then a natural choice for the
area process $\mathbb{X}^2$ is the \emph{geometric} one given by
\begin{equation}
  \label{eq:geom-lifting}
(\mathbb{X}^{2}_{\text{geom}})^{ij}_{\xi\eta} = \int_\xi^\eta
(X_\rho-X_\xi)^i dX_\rho^j  
\end{equation}
which naturally satisfy eq.~(\ref{eq:compatibility}) (as can be
directly checked) and eq.~(\ref{eq:area-reg}) (using lemma~\ref{lemma:young}).
\end{remark}

As shown by Lyons~\cite{Lyons}, when $\gamma > 1/3$ any integral of the form
$$
\int \varphi(X) dX
$$
can be defined to depend in a continuous way on the $\gamma$-rough path
$(X,\mathbb{X}^2)$ for sufficiently regular $\varphi$.

In~\cite{Gubinelli} it is pointed out that any $\gamma$-rough path $\mathbb{X}$ define  a
natural class of paths for which path-integrals are meaningful.
Define the Banach space $\DX$ of paths \emph{weakly-controlled by $X$} as the set of paths $Y$ that can be
decomposed as
\begin{equation}
\label{eq:decomposition}
Y_\xi - Y_\eta = Y'_{\eta}(X_\xi - X_\eta) + R^Y_{\eta\xi}
\end{equation}
with $Y' \in C^\gamma([0,1], \R^3 \otimes \R^3)$ and $\norm{R^Y}_{2\gamma} < \infty$.
Define the norm for  $Y \in \DX$
as
$$
\norm{Y}_{D} := \norm{Y'}_{\gamma}+ \norm{R^Y}_{2\gamma} + \norm{Y'}_\infty;
$$
moreover let
$$
\norm{Y}^*_{D} := \norm{Y}_{D}+ \norm{Y}_{\infty}
$$

Since we will need to consider only closed paths we will require for $Y
\in \DX$ that $Y_0 = Y_1$. Then it is easy to show that
$$
\|Y\|_\gamma \le \|Y\|_{D} (1+\|X\|_\gamma)
$$
and that $\DX \subseteq \CC^\gamma$.

The next lemma states that $\DX$ behaves nicely under maps by regular functions:
\begin{lemma}
\label{lemma:stability-D}
If $\varphi$ is a $C^2$ function and $Y \in \DX$ then $\varphi(Y) \in
\DX$ with $\varphi(Y)' = \nabla \varphi(Y) Y'$ and there exists a
constant $K \ge 1$ such that
\begin{equation}
  \label{eq:bound-varphi-y}
\|\varphi(Y)\|_D \le K \|\nabla \varphi\|_1 \|Y\|_D (1+\|Y\|_D) (1+\|X\|_\gamma)^2.
\end{equation}
Moreover if $Y, \wt Y \in \DX$ we have
\begin{equation}
  \label{eq:bound-varphi-y-diff}
\|\varphi(Y)-\varphi(\wt Y)\|_D \le K \|\nabla \varphi\|_2 \|Y\|_D
(1+\|Y\|_D+\|\wt Y\|_D)^2 (1+\|X\|_\gamma)^4 \|Y-\wt Y\|_D.
\end{equation}
\end{lemma}
\begin{proof}
See~\cite[Prop. 4]{Gubinelli}.  
\end{proof}

The main result about weakly-controlled paths is that they can be
integrated one against the other with a good control of the resulting
object:

\begin{lemma}[Integration of weakly-controlled paths]
\label{lemma:rough-integrals}
If $Y,Z \in \DX$ then the integral
$$
\int_\xi^\eta Y dZ  := Y_\xi (Z_\eta-Z_\xi) + Y'_\xi Z'_\xi
\mathbb{X}^2_{\eta\xi} + Q_{\xi\eta},\qquad \eta,\xi \in [0,1]
$$
is well defined
with
$$
\norm{Q}_{3\gamma} \le  C_\gamma' C_X \norm{Y}_{D} \norm{Z}_{D}
$$
where $C'_\gamma > 1$ and
$$
C_X =  (1+\norm{X}_\gamma+\norm{\mathbb{X}^2}_{2\gamma}).
$$

The integral $\int_\xi^\eta Y dZ$ is the limit of the
following ``renormalized'' finite sums
$$
\sum_{i=0}^{n-1}\left[ Y_{\xi_i} (Z_{\xi_{i+1}}- Z_{\xi_{i}}) + Y'_{\xi_i} Z'_{\xi_i} \mathbb{X}^2_{\xi_{i+1}\xi_i}\right]
$$
(where $\xi_0=\xi < \xi_1 < \cdots \xi_n = \eta$ is a finite partition
of $[\xi,\eta]$) as the size of the partition goes to zero.

Moreover if 
$\wt Y,\wt Z \in \DWTX$, then
$$
\int_\xi^\eta Y dZ-\int_\xi^\eta \wt Y d\wt Z  = Y_\xi (Z_\eta-Z_\xi)
- \wt Y_\xi (\wt Z_\eta-\wt Z_\xi) 
+( Y'_\xi Z'_\xi
-\wt Y'_\xi \wt Z'_\xi)
\mathbb{ X}^2_{\eta\xi} + Q_{\xi\eta}- \wt Q_{\xi\eta}
$$
and
$$
\norm{Q-\wt Q}_{3\gamma} \le  C'_\gamma C_X (\norm{Y}_{D} \epsilon_Y+
\norm{Z}_{D} \epsilon_Z)
$$
with
$$
\epsilon_Y = \norm{Y'-\wt Y'}_\infty +\norm{Y'-\wt Y'}_\gamma + \norm{R^Y  -R^{\wt Y}}_{2\gamma}
$$
$$
\epsilon_Z = \norm{Z'-\wt Z'}_\infty +\norm{Z'-\wt Z'}_\gamma + \norm{R^Z  -R^{\wt Z}}_{2\gamma}
$$
\end{lemma}
\begin{proof}
See~\cite[Theorem 1]{Gubinelli}.  
\end{proof}

A weakly-controlled path $(Y,Y') \in \DX$ can be naturally lifted to a
$\gamma$-rough path by setting
\begin{equation}
  \label{eq:lifting}
  \mathbb{Y}^{2,ij}_{\rho\xi} \coloneq \int_{\rho}^\xi
  (Y^i_{\eta}-Y^i_{\rho}) dY^j_{\eta}
\end{equation}
where the integral is the rough integral in Lemma~\ref{lemma:rough-integrals}.

\subsection{Local existence and uniqueness}

Given $T > 0$,
consider the Banach space $\DXT = C([0,T],\DX)$ endowed with the norm
$$
\norm{Y}_{\DXT} \coloneq \sup_{t \in [0,T]} \norm{Y(t)}_{D}^*
$$
and, as above, the application $F : \DXT \to \DXT$ defined as
$$
F(Y)(t)_\xi \coloneq X_{\xi} + \int_0^t V^{Y(s)}(Y(s)_\xi) ds
$$
with
\begin{equation}
  \label{eq:rough-V}
V^Y(x) \coloneq \int_0^1 A(x-Y_{\eta}) dY_\eta  
\end{equation}
understood as a rough integral. 
Lemma~\ref{lemma:stability-D} guarantees, under suitable
smoothness of the function $x \mapsto V^{Y(s)}(x)$, that $F(Y) \in \DXT$ if $Y \in \DXT$ with
$F(Y)'$ given by 
\begin{equation}
  \label{eq:Fderivative}
[F(Y)(t)']_\xi^{ij} \coloneq \delta_{ij} + \sum_{k=1}^3 \int_0^t \nabla_k
[V^{Y(s)}(Y(s)_\xi)]^i
[Y(s)'_\xi]^{kj}\, ds, \qquad i,j=1,2,3,
\end{equation}
where $\delta_{ij}$ is the Kronecker symbol.

A result proven in~\cite[p.~103]{Gubinelli} implies that the rough
integral
in eq.(\ref{eq:bound-rough-v})
can be indifferently understood as an integral with respect to  the
weakly-controlled path $Y \in \DX$ or as an integral over the
lifted rough-path $(Y,\mathbb{Y}^2)$ where $\mathbb{Y}^2$ is defined
as in eq.(\ref{eq:lifting}).

\medskip

We will state now the main result of this section, namely the
existence and uniqueness of solutions to the vortex line equation in
the space $\DX$.

\begin{theorem}
\label{th:sol-rough}
Assume $\|\nabla A\|_4 < \infty$ and $\mathbb{X}$ is a $\gamma$-rough
path. Then there
exists a time $T_0 >0 $ depending only on $\|\nabla A\|_3, X,\mathbb{X},
\gamma$ such that the equation~(\ref{eq:integrated}) has a unique
solution bounded in $\DX$ for any $T \le  T_0$.
\end{theorem}
\begin{proof}
 Lemma~\ref{lemma:existence-rough} and lemma~\ref{lemma:uniqueness-rough} prove that on a small enough time
interval $[0,T]$ the map $F$ is a strict contraction on some ball of $\DXT$ having a unique fixed
point. The arguments are similar to those used in the case $\gamma > 1/2$. 
\end{proof}

\medskip

Since $V$ is defined through rough integrals we can obtain the
following bounds on its regularity:

\begin{lemma}
\label{lemma:bounds-rough}
Let $Y,\wt Y \in \DX$
, for any integer $n \ge 0$:
  \begin{equation}
    \label{eq:bound-rough-v}
\norm{\nabla^n V^Y} \le
4 C'_\gamma
\norm{\nabla^{n+1} A}_1 C_X^3 \norm{Y}_D^2 (1+\norm{Y}_D)
  \end{equation}
and
  \begin{equation}
    \label{eq:bound-rough-v2}
\norm{\nabla^n V^Y-\nabla^n V^{\wt Y}} \le
16 C'_\gamma C_X^3 \norm{\nabla^{n+1} A}_2 \norm{Y-\wt Y}_D^*
  (1+\norm{Y}_D)^2 \norm{Y}_D
  \end{equation}
\end{lemma}
\begin{proof}
See Appendix, section~\ref{sec:proof-bounds-rough}.  
\end{proof}

\begin{lemma}
\label{lemma:existence-rough}
Assume $\|\nabla A\|_3 < \infty$.
For any initial $\gamma$-rough path $\mathbb{X}$ with $\gamma > 1/3$
there exists a time $T_0$ such that for any
 time $T \le T_0$ the set
$$
Q_T \coloneq \{Y \in \DXT : Y(0)=X, \norm{Y}_{\DXT} \le B_T \}
$$
where $B_T$ is a suitable constant, is invariant under $F$.
\end{lemma}
\begin{proof}
Fix a time $T > 0$.
First of all we have, for any $t \in [0,T]$
\begin{equation*}
  \begin{split}
|F(Y)(t)_\xi| &\le |X_{\xi}| + \int_0^t \norm{V^{Y(s)}}_\infty ds
    \\ & \le \norm{X}_\infty + 4 C'_\gamma
\norm{\nabla A}_1 C_X^3  \int_0^t ds \norm{Y(s)}_D^2 (1+\norm{Y(s)}_D)
\\ & \le
\norm{X}_\infty +  4 C'_\gamma T C_X^3
\norm{\nabla A}_1  \norm{Y}_{\DXT}^2 (1+\norm{Y}_{\DXT})
  \end{split}
\end{equation*}
so that
$$
\sup_{t\in[0,T]} \norm{F(Y)(t)} \le
\norm{X}_\infty + 4 T C'_\gamma  C_X^3
\norm{\nabla A}_1  \norm{Y}_{\DXT}^2 (1+\norm{Y}_{\DXT}).
$$

Next,
\begin{equation*}
  \begin{split}
\|F(Y)(t)\|_D &\le \|X\|_D + \int_0^t \norm{V^{Y(s)}(Y(s)_\cdot)}_D ds
\\ &\le \|X\|_D + K C_X^2 \int_0^t \norm{\nabla V^{Y(s)}}_1 \norm{Y(s)}_D (1+\norm{Y(s)}_D) ds
  \end{split}
\end{equation*}
where we used Lemma~\ref{lemma:stability-D}.

Lemma~\ref{lemma:bounds-rough} gives then
\begin{equation*}
  \begin{split}
\|F(Y)(t)\|_D &\le \|X\|_D + 16 K C'_\gamma C_X^5 \int_0^t  \norm{Y(s)}^3_D (1+\norm{Y(s)}_D)^2 ds
  \end{split}
\end{equation*}
from which we easily obtain
\begin{equation*}
  \begin{split}
  \norm{F(Y)}_{\DXT}
& \le \|X\|_D^* + C'_\gamma T \left [16 K  C_X^5 \norm{Y}_{\DXT}^3
  (1+\norm{Y}_{\DXT})^2 + 4  C_X^3 \norm{Y}_{\DXT}^2
  (1+\norm{Y}_{\DXT}) \right] 
\\ &\le \norm{X}_{\infty} + 1 +
20 K C'_\gamma C_X^5 T \norm{\nabla A}_3 \norm{Y}_{\DXT}^2 (1+\norm{Y}_{\DXT})^3
      \end{split}
\end{equation*}
and for $T$ small enough ($T \le T_0$ with $T_0$ depending only on
$\mathbb{X}$ and $\norm{\nabla A}_3$) we have that there exists a constant $B_T$
such that if $\norm{Y}_{\DXT} \le B_T$ then  $\norm{F(Y)}_{\DXT} \le B_T$.
\end{proof}

\begin{lemma}
\label{lemma:uniqueness-rough}
 Assume $\|\nabla A\|_4 < \infty$.
There exists a time $\overline{T} \le T_0$ such that the
map $F$ is a strict contraction on $Q_{\overline{T}}$.
\end{lemma}
\begin{proof}
Let $Z(t) := F(Y)(t)$
and $\wt Z (t) := F(\wt Y)(t)$, with $Y, \wt Y \in Q_{T}$.
Let $H := Z - \wt Z$.
We start with the estimation of the sup norm of $H(t)$:
\begin{equation*}
  \begin{split}
 H(t)_\xi & =   Z(t)_\xi-\wt Z(t)_\xi   = \int_0^t ds \left[
V^{Y(s)}(Y(s)_\xi) - V^{\wt Y(s)}(\wt Y(s)_\xi)
\right]
\\ & =
\int_0^t ds \left[
V^{Y(s)}(Y(s)_\xi)- V^{Y(s)}(\wt Y(s)_\xi) + V^{Y(s)}(\wt Y(s)_\xi) - V^{\wt Y(s)}(\wt Y(s)_\xi)
\right]
  \end{split}
\end{equation*}
which gives
\begin{equation}
\label{eq:bound-h0}
  \begin{split}
 \norm{H(t)}_\infty & \le \int_0^t ds \left [\norm{\nabla V^{Y(s)}}_\infty
 \norm{Y(s)-\wt Y(s)}_\infty  + \norm{V^{Y(s)}-V^{\wt Y(s)}}_\infty
 \right]
\\ & \le C'_\gamma T \norm{Y-\wt Y}_{\DXT} \left[4 C_X^3 \norm{\nabla^2 A}_1
 \norm{Y}_{\DXT}^2 (1+\norm{Y}_{\DXT})
 \right. \\
& \qquad  \left. + 16 C_X^3 \norm{\nabla A}_2
 \norm{Y}_{\DXT} (1+\norm{Y}_{\DXT})^2 \right]
\\ & \le 20 C'_\gamma T C_X^3 \norm{\nabla A}_2 (1+\norm{Y}_{\DXT})^3 \norm{Y-\wt Y}_{\DXT}.
  \end{split}
\end{equation}

Next we need to estimate the $\DX$ norm of $H(t)$. 

\begin{equation*}
  \begin{split}
\norm{H(t)}_{D} & \le  \int_0^T \norm{V^{Y(s)}(Y(s)_\cdot) -
V^{\wt Y(s)}(\wt Y(s)_\cdot)}_D \,ds
\\ & \le  \int_0^T \left[\norm{V^{Y(s)}(Y(s)_\cdot) - V^{\wt
    Y(s)}(Y(s)_\cdot)}_D  +  \norm{V^{\wt Y(s)}(Y(s)_\cdot)-
V^{\wt Y(s)}(\wt Y(s)_\cdot)}_D \right]ds
  \end{split}
\end{equation*}
Next we estimate the first contribution in the integral by using
Lemma~\ref{lemma:stability-D} and Lemma~\ref{lemma:bounds-rough} as
\begin{equation*}
  \begin{split}
\norm{V^{Y(s)}(Y(s)_\cdot) - & V^{\wt
    Y(s)}(Y(s)_\cdot)}_D  \le
K C_X^2 \norm{\nabla V^{Y(s)} - \nabla V^{\wt
    Y(s)} }_1 \norm{Y(s)}_D (1+\norm{Y(s)}_D)
\\ & \le 32
K C'_\gamma C_X^5 \norm{\nabla^2 A}_3 \norm{Y(s)-\wt Y(s)}^*_D \norm{Y(s)}_D^2 (1+\norm{Y(s)}_D)^3
\\ & \le 32
K C'_\gamma C_X^5 \norm{\nabla^2 A}_3 \norm{Y-\wt Y}_{\DXT} \norm{Y}_{\DXT}^2 (1+\norm{Y}_{\DXT})^3
  \end{split}
\end{equation*}
and the second as
\begin{equation*}
  \begin{split}
 \norm{V^{\wt Y(s)}(Y(s)_\cdot)- &
V^{\wt Y(s)}(\wt Y(s)_\cdot)}_D
 \le K C_X^4 \norm{\nabla V^{\wt Y(s)}}_2 (1+\norm{\wt
    Y(s)}_D+ \norm{Y(s)}_D)^2 \norm{Y-\wt Y}_D     
\\ & \le 8 K C'_\gamma \norm{\nabla^2 A}_3 C_X^7  (1+\norm{\wt
    Y(s)}_D+ \norm{Y(s)}_D)^3 \norm{\wt Y(s)}^2_D \norm{Y(s)-\wt Y(s)}_D     
\\ & \le 8 K C'_\gamma \norm{\nabla^2 A}_3 C_X^7  (1+\norm{\wt
    Y}_{\DXT}+ \norm{Y}_{\DXT})^3 \norm{\wt Y}^2_{\DXT} \norm{Y-\wt Y}_{\DXT}     
  \end{split}
\end{equation*}
giving
\begin{equation*}
  \begin{split}
\norm{H(t)}_{D} &  
\le  40 T  K C'_\gamma \norm{\nabla^2 A}_3 C_X^7 
(1+\norm{Y}_{\DXT}+\norm{Y}_{\DXT})^5 \norm{Y-\wt Y}_{\DXT}
  \end{split}
\end{equation*}

Finally, collecting together the
bounds

\begin{equation*}
  \norm{F(Y)-F(\wt Y)}_{\DXT}  
\le 60 T C'_\gamma K C_X^7 \norm{\nabla A}_4 (1+\norm{Y}_{\DXT}+\norm{\wt Y}_{\DXT})^5
\norm{Y-\wt Y}_{\DXT}  
\end{equation*}

When $Y, \wt Y \in Q_T$ we have
\begin{equation*}
  \norm{F(Y)-F(\wt Y)}_{\DXT}  
\le 60 T C'_\gamma K C_X^7 \norm{\nabla A}_4 (1+2 B_T)^5
\norm{Y-\wt Y}_{\DXT}  
\end{equation*}
Proving that for $\overline{T} < T_0$ small enough so that
\begin{equation*}
60 \overline{T} C'_\gamma K C_X^7 \norm{\nabla A}_4 (1+2 B_{\overline T})^5
< 1
\end{equation*}
$F$ is a contraction in the ball $Q_{\overline{T}}\subset\DXoT$.
\end{proof}

\begin{remark}
By imposing enough regularity on $A$ and requiring that $X$
can be completed to a geometric rough path with bounded $p$-variation
(in the sense of Lyons)  it is likely that the above proof of
existence and uniqueness can be extended to cover the case of rougher initial conditions
(e.g. paths living in some $\CC^\gamma$ with $\gamma < 1/3$ suitably
lifted to rough-paths).
\end{remark}

\begin{remark}
The solution $(Y,Y')$ in $\mathbf{D}_{X,\overline{T}}$, satisfy
 (compare with eq.~(\ref{eq:Fderivative}))
\begin{equation}
  \label{eq:sol1}
Y(t)^{\prime}_\eta = \id + \int_0^t  \nabla  V^{Y(s)}(Y(s)_{\eta}) Y^{\prime}(s)_\eta  \,ds
\end{equation}
and, as can be easily verified,
\begin{equation}
  \label{eq:sol2}
  \begin{split}
[R^{Y}(t)]^i_{\xi\eta} & = \int_0^t ds \left[\Big.
\nabla^k
V^{Y(s)\,i}(Y(s)_{\eta}) [R^{Y}(s)]^k_{\xi\eta}     +\right.
\\ & \quad + \left.\int_0^1 dr \int_0^r dw \sum_{k,l = 1,2,3}\nabla^k \nabla^l
  [V^{Y(s)}(Y(s)_\xi + w Y(s)_{\xi\eta} )]^i
Y(s)^k_{\xi\eta} Y(s)^l_{\xi\eta}      \right].
  \end{split}
\end{equation}  
where $Y(s)_{\xi\eta} \coloneq Y(s)_\eta - Y(s)_\xi$.
\end{remark}


\subsection{Dependence on the initial condition}
Let $\Sigma^\gamma$ the set of $\gamma$-rough paths $(X,\mathbb{X}^2)$
(with $\gamma > 1/3$). On $\Sigma^\gamma$ consider the distance
$$
d(\mathbb{X},\mathbb{\wt X}) = \|X-\wt X\|_\gamma +
\|\mathbb{X}^2-\mathbb{\wt X}^2\|_{2\gamma}
$$
where $\mathbb{X} = (X,\mathbb{X}^2)$ and $\mathbb{\wt X} = (\wt
X,\mathbb{\wt X}^2)$ are two points in $\Sigma^\gamma$.

Fix $r>0$ and let $\mathscr{B}(\Sigma^\gamma;r)$ the open ball of
$\Sigma^\gamma$ with radius $r$ centered at the trivial rough-path
$(0,0)$. Let $T_0(r)>0$ the smaller existence time of solution to the evolution
problem, as given by Theorem~\ref{lemma:existence-rough}, for initial
data living in $\mathscr{B}(\Sigma^\gamma;r)$. 

Define the map $\Gamma: \mathscr{B}(\Sigma^\gamma;r) \to
C([0,T_0],\Sigma^\gamma)$ as follows: for any $t \in [0,T_0]$ and
any $\mathbb{X} = (X,\mathbb{X}^2) \in \mathscr{B}(\Sigma^\gamma;r)$ let
$\Gamma(\mathbb{X}) \coloneq \mathbb{Y}$ where
$\mathbb{Y} = (Y,\mathbb{Y}^2)$ with $Y\in \DXTO$ the solution
starting at $\mathbb{X}$ and $\mathbb{Y}^2$ the corresponding area
process defined as 
\begin{equation}
  \label{eq:lifting2}
  [\mathbb{Y}^{2}(t)]^{ij}_{\rho\xi} \coloneq \int_{\rho}^\xi
  (Y(t)^i_{\eta}-Y(t)^i_{\rho}) dY(t)^j_{\eta}
\end{equation}
for any $t \in [0,T_0]$ where the integral is understood as an
integral over the weakly-controlled path $Y(t) \in \DX$.

Then
\begin{theorem}
The map $\Gamma$ is Lipshitz from $\mathscr{B}(\Sigma^\gamma;r)$
to
$C([0,T_0],\Sigma^\gamma)$ endowed with the uniform distance.
\end{theorem}
\begin{proof}
Take two inital conditions in $\mathscr{B}(\Sigma^\gamma;r)$:
$\mathbb{X} = (X,\mathbb{X}^2)$ and $\mathbb{\wt X} = (\wt
X,\mathbb{\wt X}^2)$. Let $Y \in \DXTO$ (resp. $\wt Y \in \DWTXTO$) the solution starting from
$\mathbb{X}$ ($\mathbb{\wt X}$).

Let 
$$
\Psi(t) \coloneq \|Y'(t)-\wt Y'(t)\|_{\gamma}^* +\|R^Y(t)-R^{\wt Y}(t)\|_{2\gamma}.
$$

Using results form~\cite{Gubinelli} it is not too difficult
to prove that (cfr. Lemma~\ref{lemma:rough-integrals} and Lemma~\ref{lemma:bounds-rough})
\begin{equation}
  \label{eq:bound-area-diff}
\|\mathbb{Y}^2(t)-\mathbb{\wt Y}^2(t)\|_{2\gamma} \le D_1 [d(X,\wt
X)+ \Psi(t)]
\end{equation}
and 
\begin{equation}
  \label{eq:bound-v-diff}
\|\nabla^n V^{Y(t)}- \nabla^n V^{\wt Y(t)}\| \le D_2 \|\nabla^{n+1} A\|_2 [d(X,\wt
X)+\Psi(t)], \qquad n \ge 0
\end{equation}
uniformly for $t \in [0,T_0]$, where here and in the following $D_k > 0$  are constants depending only
on $r$ and on $\gamma$.

At this point $\Psi(t)$ can be estimated using the expression given in
eq.(\ref{eq:sol1}) and eq.(\ref{eq:sol2}) for $Y'(t),\wt Y'(t),
R^{Y}(t)$ and $R^{\wt Y}(t)$ to obtain
\begin{equation*}
\Psi(t) \le  D_3 \|\nabla A\|_4\int_0^t  [d(\mathbb{X},\mathbb{\wt
  X})+ \Psi(s)]   ds
\end{equation*}
For $\overline T$ small enough so that
$
 D_3 \|\nabla A\|_4 \overline{T} \le 1/2
$
we have
$$
\sup_{t \in [0,\overline{T}]}\Psi(t) \le d(\mathbb{X},\mathbb{\wt
  X})
$$
which implies that
$$
 \sup_{t \in [0,\overline{T}]} d(\mathbb{Y}(t),\mathbb{\wt Y}(t)) \le
 D_4 d(\mathbb{X},\mathbb{\wt
  X})
$$
(using eq.(\ref{eq:bound-area-diff})).
Then again a simple induction argument allows to extend this result to
the whole interval $[0,T_0]$ proving the claim.
\end{proof}

The continuity of $\Gamma$ implies in particular that if we have a
sequence of smooth loops $(X^{(n)})_{n\ge 1}$ which can be naturally
lifted to a sequence of $\gamma$-rough paths $(\mathbb{X}^{(n)})_{n
  \ge 1}$ (using forumla~(\ref{eq:geom-lifting})
for the area process $\mathbb{X}^{(n),2}$) and such that it converges to
the rough path $\mathbb{X}$ (in the topology of $\Sigma^\gamma$) then the sequence of
solutions $\Gamma(\mathbb{X}^{(n)})$ converge 
 to the solution  $\Gamma(\mathbb{X})$.

Since  integrals over a smooth path $X^{(n)}$ coincide with integrals over
the (geometrically lifted) rough path $\mathbb{X}^{(n)}$ we have that for smooth initial
conditions classical solutions to the vortex line equation coincide
with the projection of the solution $\Gamma(\mathbb{X}^{(n)})$ built in the space of rough-paths.   

Then the sequence of classical solutions (once lifted to a $\gamma$-rough-path)  converges
in the sense of $\gamma$-rough paths (which is stronger than the
$\gamma$-H\"older topology) to the solution  $\Gamma(\mathbb{X})$.

\subsection{Dynamics of the covariations}
\label{sec:structure}
Recall the framework described in Sec.~\ref{sec:quad-var} on the
covariation structure of the solution. If we assume that the initial
condition $(X,\mathbb{X}^2)$ is a random variable a.s. with values in
the space of $\gamma$-rough paths (with $\gamma > 1/3$) and that $X$ is a
process with all its mutual covariations, then the solution $Y(t)$ at
any instant of time $t$ less that a random time $T_0$
(depending on the initial condition)  it is still a process with all
its mutual covariations (due to  Prop.\ref{PQVar}).    

The covariations of $Y$ satisfy the equation
\begin{equation}
  \label{eq:Yqvar}
  [ Y(t)^i,Y(t)^j ]_{\eta} = \sum_{k,l=1,2,3}
\int_{0}^{\eta} (Y(t)')^{ik}_\rho (Y(t)')^{jl}_\rho d_{\rho} [ X^k, X^l]_{\rho}
\end{equation}
Indeed, comparing
eq.~(\ref{eq:sol1}) with eq.~(\ref{eq:Mprocess}) we can identify the
function $M(t)_\xi$ in Prop.~\ref{PQVar} with $Y(t)'_\xi$.
\begin{remark}
The same result can be obtained noting that, for our solution,
\begin{equation*}
 \sum_i |Y(t)_{\xi_{i+1} \xi_{i}}|^2 = \sum_i Y(t)'_{\xi_i} X_{\xi_{i+1}
 \xi_{i}} Y'(t)_{\xi_i} X_{\xi_{i+1} \xi_{i}} + \sum_i O(|\xi_{i+1}-\xi_i|^{3\gamma})
.
\end{equation*}  
\end{remark}

Eq.~(\ref{eq:Yqvar}) has a differential counterpart in the following equation
\begin{equation}
  \label{eq:dyn-Yqvar}
\frac{d}{dt} W(t)_\xi = \int_0^\xi (H(t)_\rho d_\rho W(t)_\rho 
+  d_\rho W(t)_\rho H(t)_\rho^*)
\end{equation}
where we let $W(t)_\xi := [Y(t),Y(t)^*]_\xi$ as a matrix valued process, and 
$H(t)_\xi := \nabla V^{Y(t)}(Y(t)_\xi)$. 
To understand better this evolution equation let us split the matrix
$H(t)_\xi$ into its symmetric $S$ and anti-symmetric $T$ components:
$$
H(t)_\xi = S(t)_\xi + T(t)_\xi, \qquad S(t)_\xi = S(t)_\xi^*,\qquad T(t)_\xi
= - T(t)_\xi^*.
$$

Moreover define $Q(t)_\xi$ as the solution of the Cauchy problem
$$
\frac{d}{dt} Q(t)_\xi = - Q(t)_\xi T(t)_\xi, \qquad Q(0)_\xi = \text{Id}
$$
i.e. 
$$
Q(t)_\xi = \exp\left[-\int_0^t T(s)_\xi ds \right].
$$
Since $T$ is antisymmetric, the matrix $Q$ is orthogonal,
i.e. $Q(t)_\xi^{-1} = Q(t)_\xi^*$. This matrix describe the
rotation of the local frame of reference at the point $Y(t)_\xi$
caused by the motion of the curve. 

Then define 
$$
\widetilde W(t)_{\xi} := \int_0^\xi Q(t)_\rho d_\rho W(t)_\rho Q(t)_\rho^{-1}
$$ 
and analogously $\widetilde
T(t)_\xi := Q(t)_\xi T(t)_\xi Q(t)_\xi^{-1}$ and $\widetilde S(t)_\xi
=  Q(t)_\xi S(t)_\xi Q(t)_\xi^{-1}$, and compute the following time-derivative:
\begin{equation*}
  \begin{split}
 \frac{d}{dt} d_\xi \widetilde W(t)_\xi & = \frac{dQ(t)_\xi}{dt} Q(t)_\xi^{-1}
 d_\xi \widetilde W(t)_\xi + d_\xi \widetilde W(t)_\xi Q(t)_\xi \frac{dQ(t)_\xi^{-1}}{dt} +
 Q(t)_\xi \left(\frac{d}{dt} d_\xi W(t)_\xi \right) Q(t)_\xi^{-1}  
\\ & = - \widetilde T(t)_\xi 
 d_\xi\widetilde W(t)_\xi + d_\xi\widetilde W(t)_\xi \widetilde T(t)_\xi+
 Q(t)_\xi [H(t)_\xi d_\xi W(t)_\xi +  d_\xi W(t)_\xi H(t)_\xi^* ] Q(t)_\xi^{-1}   
\\ & = \widetilde S(t)_\xi d_\xi \widetilde W(t)_\xi +
 d_\xi \widetilde W(t)_\xi \widetilde S(t)_\xi
  \end{split}
\end{equation*}
This result implies that $d_\xi W(t)_\xi$ can be decomposed as
\begin{equation}
  \label{eq:stretching}
d_\xi W(t)_\xi = Q(t)_\xi^{-1} \exp\left[\int_0^t \widetilde S(s)_\xi ds\right]
d_\xi W(0)_\xi
\exp\left[\int_0^t \widetilde S(s)_\xi ds\right]^* Q(t)_\xi  
\end{equation}
The relevance of this decomposition is the following. Modulo
rotations, the matrix $\widetilde S(t)_\xi$, corresponds to the
symmetric part of the tensor field $\nabla V^{Y(t)}(x)$ in the point
$x = Y(t)_\xi$. This symmetric component describe the stretching of
the  volume element around $x$ due to the flow generated by
the (time-dependent)  vector
field $V^{Y(t)}$. The magnitude of the covariation then varies
with time, due to this
 stretching contribution, according to eq.~(\ref{eq:stretching}).
   


\section{Random vortex filaments}
\label{sec:random}

\subsection{Fractional Brownian loops with $H > 1/2$}
\label{sec:fbl}

Consider the following probabilistic model of Gaussian vortex
filament. Let $(\wt X_\xi)_{\xi \in [0,1]}$ a 3d fractional
Brownian Motion (FBM) of Hurst index $H$, i.e. a centered Gaussian process on
$\mathbb{R}^3$ defined on the probability space
$(\Omega,\mathbb{P},\mathcal{F})$ such that 
$$
 \expect \wt X^i_\xi \wt X^j_\eta 
= \frac{\delta_{ij}}{2}  (|\xi|^{2H}+|\eta|^{2H}-|\xi-\eta|^{2H}),
 \qquad
i,j=1,2,3,\; \xi,\eta \in[0,1]
$$
with $H > 1/2$ and $X_0 = 0$. Define the Gaussian process $(X_\xi)_{\xi\in[0,1]}$
 as 
\begin{equation}
  \label{eq:fbm-bridge}
X_\xi \coloneq \wt X_\xi - \frac{C(\xi,1)}{C(1,1)} \wt X_1  
\end{equation}
where $C(\xi,\eta) \coloneq (|\xi|^{2H}+|\eta|^{2H}-|\xi-\eta|^{2H})$.
Then $X_0 = X_1 = 0$ a.s. moreover the process $(X_\xi)_\xi$ is independent of
the r.v. $\wt X_1$. We call $X$ a fractional Brownian loop (FBL).
Using the standard Kolmogorov criterion it is easy to show that $X$ in
a.s. H\"older continuous for any index $\gamma < H$. Since $H > 1/2$
then we can choose $\gamma \in(1/2,H)$ and apply the results of
Sec.~\ref{sec:young} to obtain the evolution of a random vortex
filament modeled on a FBL.

\subsection{Evolution of Brownian loops}
\label{sec:brownian}
As an example of application of Theorem~\ref{th:sol-rough} we can consider
the evolution of an initial random curve whose law is that of a
Brownian Bridge on $[0,1]$ starting at an arbitrary point $x_0$. A standard
three-dimensional Brownian Bridge $\{B_\xi\}_{\xi \in[0,1]}$ such that $B_0 =
B_1 = x_0 \in \mathbb{R}^3$ is a stochastic process defined on a
complete probability space $(\Omega,\mathcal{F},\prob)$ whose law is the law of
a Brownian motion starting at $x_0$ and conditioned to reach $x_0$ at
``time'' $1$. As in the previous section, it can be obtained starting
from a standard Brownian motion $\{ \wt B\}_{\xi \in [0,1]}$ as
$$
B_\xi = \wt B_{\xi} - \xi \wt B_1, \qquad \xi \in [0,1].
$$

The Brownian Bridge is a semi-martingale with respect to its
own filtration $\{\mathcal{F}^B_\xi : 0 \le \xi \le 1\}$ with
decomposition
$$
dB_\xi = \frac{B_\xi - x_0}{1-\xi} d\xi + d\beta_\xi
$$
where $\{\beta_\xi\}_{\xi \in [0,1]}$ is a standard 3d Brownian motion. Using
the results in~\cite{Gubinelli}, it is easy to see that $B$ is a
$\gamma$-H\"older rough path if we consider it together with the area
process defined as
\begin{equation}
  \label{eq:loop-area}
\mathbb{B}^{2,ij}_{\xi\eta} = \int_\xi^\eta (B^i_\rho-B^i_\eta) \circ dB_\rho^j  \end{equation}
where the integral is understood in Stratonovich sense. Indeed there
exists a version of the process
$(\xi,\eta) \mapsto \mathbb{B}^2_{\xi\eta}$ which is continuous in
both parameters and such that $\|\mathbb{B}\|_{2\gamma}$ is almost
surely finite (also all moments are finite). Then outside an event of
$\prob$-measure zero the couple $(B,\mathbb{B})$ is a
$\gamma$-H\"older rough path and by theorem~\ref{th:sol-rough} there
exists a solution of the problem~(\ref{eq:evolution}) starting at
$B$. Of course, in this case, the solution depends a priori on the
choice~(\ref{eq:loop-area}) we made for the area process. Indeed if in
eq.~(\ref{eq:loop-area}) we consider, for example, the It\^o integral
(for which the regularity result on $\mathbb{B}$ still holds) we would
have obtained a different solution, even if the path $B$ is unchanged.

\begin{remark}
Consider the discussion in Sec.~\ref{sec:structure} and note that,
for our Brownian loop $B$ the covariation is
$[B,B]_\xi =  \text{Id}\cdot \xi$ we can say
that the covariation of the solution $Y$ starting at $B$ will
be
$$
d_\xi [Y(t),Y(t)^*]_\xi = 
 Q(t)_\xi^{-1} \exp\left[\int_0^t \widetilde S(s)_\xi ds\right]
\exp\left[\int_0^t \widetilde S(s)_\xi ds\right]^* Q(t)_\xi 
$$
(the notations are the same as in Sec.~\ref{sec:structure}).  
\end{remark}


\subsection{Evolution of fractional Brownian loops with $1/3 < H \le 1/2$}
\label{sec:flb2}
The above results on the Browian loop are a particular case of the more
general  case of a fractional Brownian loop $X$ of Hurst index $H \in
(1/3,1/2]$.
In general lifting $X$ to a $\gamma$-rough path (with $1/3 < \gamma
< H$) require to build an area process $\mathbb{X}^2$ with appropriate
regularity (when $H \neq 1/2$ cannot be obtained by semi-martingale
stochastic calculus as in section~\ref{sec:brownian} above). 

In~\cite{coutin} the authors give  a construction of
the \emph{area process} $ \mathbb{\wt X}^2$ in the case where
$\wt X$ is a FBM of Hurst index $H > 1/4$. Moreover they prove that
the sequence $(\wt X^{(n)}, \mathbb{\wt X}^{(n),2})_{n \in \mathbb{N}}$ of piece-wise linear
approximations of $\wt X$ toghether with the associated geometric area
process $\mathbb{\wt X}^{(n),2}$ converges to $(\wt X, \mathbb{\wt X}^{2})$ in
the generalized $p$-variation sense for any $1/H < p < 4 $ (for the
definition of this kind of convergence, see~\cite{coutin}). It is 
not difficult to prove that we have convergence also as
$\gamma$-rough-paths for any $1/3 < \gamma < H$, i.e. that
$$
\|\wt X^{(n)}-\wt X\|_\gamma + \| \mathbb{\wt X}^{(n),2} -  \mathbb{\wt X}^{2}\|_{2\gamma}
\to 0
$$
as $n \to \infty$.

To identify an appropriate area process for the fractional Brownian
loop $X$ we can consider the following definition
\begin{equation}
  \label{eq:new-area}
\mathbb{X}^2_{\xi\rho} \coloneq \mathbb{\wt X}^2_{\xi\rho} 
+ \int_\xi^\rho
(h_\eta-h_\xi) \otimes d\wt X_\eta 
+ \int_\xi^\rho
(\wt X_\eta-\wt X_\xi) \otimes dh_\eta 
+  \int_\xi^\rho
(h_\eta-h_\xi) \otimes dh_\eta   
\end{equation}
where
$$
h_\xi \coloneq - \frac{C(\xi,1)}{C(1,1)} \wt X_1
$$
(cfr. eq.(\ref{eq:fbm-bridge})). The function $h$ is $2H$-H\"older
continuous so the integrals in eq.(\ref{eq:loop-area}) can be
understood as Young integrals when $2H + \gamma > 3 \gamma > 1$.

It is then straighforward to check that $\mathbb{X}^2$ satisfy
 equation~(\ref{eq:compatibility}) and that $(X,\mathbb{X}^2)$ is a
 $\gamma$-rough path for any $1/3 < \gamma < H$.
 
Moreover by exploiting the continuity of the Young integral and the
results in~\cite{coutin} mentioned above we have that piece-wise
linear approximations $(X^{(n)},\mathbb{X}^{(n),2})$ of $(X,\mathbb{X}^2)$
converge to $(X,\mathbb{X}^2)$ as $\gamma$-rough paths.

\begin{remark}
This construction of the area process of a fractional Brownian loop
with $H > 1/3$ is a particular case of a more general result about
translations on the space of rough paths~\cite{LyonsBook,feyel,fv}.
\end{remark}

\section*{acknowledgement}
The authors would like to thank the anonymous referee for her comments
and suggestions which helped to substantially improve the paper.


\appendix
\section{Proofs of some lemmas}
\label{sec:app-proofs}

In the proofs we will need often to use Taylor expansions with integral
remainders, so for convenience we introduce a special notation:
given $X \in \CC$ let $X_{\eta \xi} \coloneq
X_\eta -X_\xi$ and $X^r_{\eta\xi } \coloneq J_r(X_\eta,X_\xi)$ where
$J_r(x,y)$ is  the linear interpolation
$$
J_r(x,y) = (x-y)r + y
$$
for $r \in [0,1]$.

\subsection{Proof of lemma~\ref{lemma:stability-young}}
\label{sec:proof-bounds-stab-young}
\begin{proof}
The bound in eq.~(\ref{eq:young-varphi-y-1}) is easy. Let us prove the
second by considering the following decomposition:
\begin{equation*}
  \begin{split}
[ \varphi(Y_\xi) - &  \varphi(\wt Y_\xi)] -  [ \varphi(Y_\eta) -
\varphi(\wt Y_\eta)]  = Y_{\xi\eta} \int_0^1 \nabla  \varphi( Y^r_{\xi\eta}) dr
 -\wt Y_{\xi\eta} \int_0^1 \nabla  \varphi(\wt Y^{ r}_{\xi\eta}) dr
\\ & = (Y_{\xi\eta} - \wt Y_{\xi\eta}) \int_0^1 \nabla  \varphi(\wt Y^r_{\xi\eta}) dr +
 Y_{\xi\eta} \left[\int_0^1 \nabla  \varphi(Y^r_{\xi\eta}) dr- \int_0^1
 \nabla \varphi(\wt Y^{
 r}_{\xi\eta}) dr\right]
\\ & = (Y_{\xi\eta} - \wt Y_{\xi\eta}) \int_0^1 \nabla  \varphi(\wt Y^r_{\xi\eta}) dr -
 Y_{\xi\eta}  \int_0^1 dr (\wt Y^{
 r}_{\xi\eta}- Y^r_{\xi\eta}) \int_0^1 dw \nabla^2  \varphi(J_w(\wt Y^{
 r}_{\xi\eta},Y^r_{\xi\eta}))
  \end{split}
\end{equation*}
where $Y^r_{\xi\eta} \coloneq J_r(Y_\xi,Y_\eta)$,
we obtain
\begin{equation*}
  \| \varphi(Y_\cdot)- \varphi(\wt Y_\cdot)\|_{\gamma} \le \|Y-\wt Y\|_{\gamma}
  \|\nabla  \varphi\| +
 \| Y\|_{\gamma} \|Y-\wt Y\|_{\infty} \|\nabla^2 A\|
\end{equation*}
which implies
\begin{equation}
  \| \varphi(Y_\cdot)- \varphi(\wt Y_\cdot)\|_{\gamma} \le 
  \|\nabla  \varphi\|_1 (1+
 \| Y\|_{\gamma})  \|Y - \wt Y\|_{\gamma}.
\end{equation}
\end{proof}

\subsection{Proof of lemma~\ref{lemma:boundsyoung}}
\label{sec:proof-bounds-young}
\begin{proof}
By eq.(\ref{eq:young-varphi-y-1}) and eq.(\ref{eq:young-varphi-y-2}) we have
\begin{equation}
\label{eq:bounda}
  \|\nabla^n A(Y_\cdot)\|_{\gamma} \le  \|\nabla^{n+1} A\| \|Y\|_{\gamma}
\end{equation}
and
\begin{equation}
\label{eq:boundadiff}
  \|\nabla^n A(Y_\cdot)-\nabla^n A(\wt Y_\cdot)\|_{\gamma} \le
  \|\nabla^{n+1} A\|_1 (1+
 \| Y\|_{\gamma})   \|Y - \wt Y\|_{\gamma}.
\end{equation}

Now, it is enough to consider $n=0$, the proof for general $n$ being
similar.
Using the lemma~\ref{lemma:young} and the bounds~(\ref{eq:bounda})
and~(\ref{eq:boundadiff}), the estimates~(\ref{eq:boundv})
and~(\ref{eq:boundvdiff}) on the velocity vector-field follow as:
\begin{equation*}
  \begin{split}
| V^{Y}(x)| = \left| \int_0^1  A(x-Y_\eta) dY_\eta \right|
   & \le  C_\gamma
 \| A(x-Y_{\cdot})\|_{\gamma} \|Y\|_{\gamma}
\\ & \le    C_\gamma  \|\nabla A\| \|Y\|_{\gamma}^2.
  \end{split}
\end{equation*}
Moreover for the difference $V^Y - V^{\wt Y}$ we have the decomposition
\begin{equation*}
  \begin{split}
    V^{Y}(x) - V^{\wt Y}(x)  & = \int_0^1 \left[
    A(x-Y_{\eta})dY_\eta - A(x-\wt Y_{\eta})d\wt
    Y_\eta\right]
\\ & =  \int_0^1
    A(x-Y_{\eta}) d (Y-\wt Y)_\eta
+ \int_0^1
    \left[ A(x-Y_{\eta})- A(x-\wt Y_{\eta})\right]
    d\wt Y_\eta
  \end{split}
\end{equation*}
which in turn can be estimated as
\begin{equation*}
  \begin{split}
| V^{Y}(x)- V^{\wt Y}(x)| & \le C_\gamma \| A(x-Y_\cdot)
\|_{\gamma} \|Y-\wt Y\|_{\gamma}
\\ & \qquad
+ C_\gamma \| A(x- Y_\cdot)-
A(x-\wt Y_\cdot) \|_{\gamma} \|\wt Y\|_{\gamma}
\\ & \le C_\gamma  \|\nabla A\| \|Y\|_{\gamma}
\|Y-\wt Y\|_{\gamma}
\\ & \qquad
+ C_\gamma  \|\wt Y\|_{\gamma} \norm{\nabla A}_1 \norm{Y-\wt
  Y}_\gamma^* (1+\norm{Y}_\gamma)
\\ & \le  C_\gamma 
\norm{\nabla A}_{1}(\norm{Y}_{\gamma}+\norm{\wt Y}_{\gamma}+\norm{\wt
  Y}_{\gamma} \norm{Y}_{\gamma}) \norm{Y-\wt Y}_{\gamma}^*
  \end{split}
\end{equation*}
giving eq.(\ref{eq:boundvdiff}).
\end{proof}

\subsection{Proof of lemma~\ref{lemma:bounds-rough}}
\label{sec:proof-bounds-rough}
\begin{proof}
Consider the case $n=0$, the general case being similar.
The path $Z_\xi = A(x-Y_\xi)$ belongs to $\DX$ and has the following
decomposition
\begin{equation*}
  \begin{split}
  Z_{\xi\eta}
& =  \nabla A(x-Y_\eta) Y_{\xi\eta} + Y_{\xi\eta} Y_{\xi\eta}
\int_0^1 dr \int_0^r dw \nabla^2 A(x-Y_{\xi\eta}^w)
\\& =  \nabla A(x-Y_\eta) Y'_\eta X_{\xi\eta}+ \nabla A(x-Y_\eta) R^Y_{\xi\eta} + Y_{\xi\eta} Y_{\xi\eta}
\int_0^1 dr \int_0^r dw \nabla^2 A(x-Y_{\xi\eta}^w)
\\ & = Z'_\eta X_{\xi\eta} + R^Z_{\xi\eta}
  \end{split}
 \end{equation*}
with
$$
Z'_\eta = \nabla A(x-Y_\eta) Y'_\eta
$$
and
$$
R^{Z}_{\xi\eta} = \nabla A(x-Y_\eta) R^Y_{\xi\eta} + Y_{\xi\eta} Y_{\xi\eta}
\int_0^1 dr \int_0^r dw \nabla^2 A(x-Y_{\xi\eta}^w)
$$
The $\DX$ norm of $Z$ can be estimated as follows:
\begin{equation}
\label{eq:boundzzz}
  \begin{split}
\norm{Z}_D & = \norm{Z'}_\infty + \norm{Z'}_\gamma +
\norm{R^Z}_{2\gamma}
\\ & \le \norm{\nabla A}_\infty \norm{Y'}_\infty + \norm{\nabla^2
  A}_\infty \norm{Y}_\gamma + \norm{\nabla A}_\infty \norm{Y'}_\gamma
\\ & \qquad + \norm{\nabla A}_\infty \norm{R^Y}_{2\gamma}+
\norm{Y}_\gamma^2 \norm{\nabla^2 A}_\infty
\\ & \le \norm{\nabla A}_1 (\norm{Y}_D + \norm{Y}_\gamma +
\norm{Y}^2_\gamma)
\\ & \le \norm{\nabla A}_1 [(1+C_X) \norm{Y}_D  +
C_X^2 \norm{Y}^2_D]
\\ & \le  C_X^2 \norm{\nabla A}_1 [2 \norm{Y}_D  +
 \norm{Y}^2_D]
  \end{split}
\end{equation}
where we used the fact that
\begin{equation}
  \label{eq:bound-gamma-D}
  \begin{split}
  \norm{Y}_\gamma & \le \norm{Y'}_\infty \norm{X}_\gamma +
  \norm{R^Y}_\gamma
\\ & \le   \norm{Y'}_\infty \norm{X}_\gamma +
  \norm{R^Y}_{2\gamma}
\\ & \le   (1+\norm{X}_\gamma) \norm{Y}_D \le C_X \norm{Y}_D .
  \end{split}
\end{equation}
Then
\begin{equation*}
V^Y(x) = \int_0^1 A(x-Y_\eta) dY_\eta = \int_0^1 Z_\eta dY_\eta
 Z_0 (Y_1-Y_0) + Z'_0 Y'_0 \mathbb{X}^2_{01} + Q_{01}
\end{equation*}
with
$$
\norm{Q}_{3\gamma} \le C'_\gamma C_X \norm{Z}_{D} \norm{Y}_D
$$
and
\begin{equation*}
  \begin{split}
  |V^Y(x)| & \le  \norm{Z'}_\infty \norm{Y'}_\infty
   \norm{\mathbb{X}^2}_{2\gamma} + \norm{Q}_{3\gamma}
\\ & \le 2 C'_\gamma C_X \norm{Z}_D \norm{Y}_D
\\ & \le 4 C'_\gamma
\norm{\nabla A}_1 C_X^3 \norm{Y}_D^2 (1+\norm{Y}_D)
  \end{split}
\end{equation*}
where we used the fact that $C'_\gamma \ge 1$.

To bound $V^Y(x)-V^{\wt Y}(x)$ we need the $\DX$ norm of the
difference $A(x-Y_\cdot)-A(x-\wt Y_\cdot)$. Let $\phi(y) = A(x-y)$
and consider the expansion
\begin{equation*}
  \begin{split}
  & \phi(Y_\eta) - \phi(Y_\xi)  - (\phi(\wt Y_\eta) - \phi(\wt Y_\xi))
\\ & =\left[ \nabla\phi(Y_{\xi})
  Y_{\eta\xi}-\nabla\phi(\wt Y_{\xi}) \wt Y_{\eta\xi}\right]
+ \int_0^1 dr\int_0^r dw \left[ \nabla^2\phi(Y_{\eta\xi}^w)
  Y_{\eta\xi} Y_{\eta\xi}-\nabla^2\phi(\wt Y_{\eta\xi}^w) \wt Y_{\eta\xi}\wt Y_{\eta\xi}\right]
  \end{split}
\end{equation*}
which by arguments similar to those leading to eq.~(\ref{eq:boundzzz})
gives a related estimate:
\begin{equation*}
  \begin{split}
 \norm{\phi(Y_\cdot)-\phi(\wt Y_\cdot)}_D \le &
\norm{\nabla \phi}_\infty \norm{Y-\wt Y}_D + \norm{\nabla^2
  \phi}_\infty \norm{Y-\wt Y}_\infty \norm{Y}_D
\\
& + 3 \norm{\nabla^3 \phi} \norm{Y-\wt Y}_\infty \norm{Y}^2_\gamma + 2
\norm{\nabla^2 \phi}_\infty  \norm{Y-\wt Y}_\gamma \norm{ Y}_\gamma
\\ & \le 6 \norm{\nabla \phi}_2 C_X^2 (1+\norm{Y}_D)^2\norm{Y-\wt Y}_D^*
  \end{split}
\end{equation*}
so that
\begin{equation*}
\norm{A(x-Y_\cdot)-A(x-\wt Y_\cdot)}_D \le 6
\norm{\nabla A}_2 C_X^2 (1+\norm{Y}_D)^2\norm{Y-\wt Y}_D^*
\end{equation*}
Now,
\begin{equation*}
  \begin{split}
V^Y(x)-V^{\wt Y}(x) & = \int_0^1 A(x-Y_\eta) dY_\eta- \int_0^1 A(x-\wt
 Y_\eta) d\wt Y_\eta
\\ & =\int_0^1 [A(x-Y_\eta)-A(x-\wt Y_\eta)]dY_\eta + \int_0^1 A(x-\wt
 Y_\eta) d(Y-\wt Y)_\eta
  \end{split}
\end{equation*}
and we can conclude by observing that
\begin{equation*}
  \begin{split}
|V^Y(x)-V^{\wt Y}(x)| & \le 2 C'_\gamma C_X
(\norm{A(x-Y_\cdot)-A(x-\wt Y_\cdot)}_D\norm{Y}_D
\\ & \qquad +\norm{A(x-Y_\cdot)}_D \norm{Y-\wt
 Y}_D)
\\ & \le  16 C'_\gamma C_X^3 \norm{\nabla A}_2 \norm{Y-\wt Y}_D^*
  (1+\norm{Y}_D)^2 \norm{Y}_D
  \end{split}
\end{equation*}

\end{proof}

\end{document}